%% file: Coordinate_Descent_Algorithms.tex
\documentclass[titlepage]{article}
\RequirePackage[colorlinks=true,linkcolor=blue,citecolor=blue,urlcolor=blue,pdfborder={0 0 0},bookmarksnumbered,pageanchor=false,plainpages=false,pdfpagelabels]{hyperref}
\usepackage{amsthm}
\usepackage{amssymb,amsmath}
\usepackage{algorithm}
\usepackage{graphicx}
\usepackage{listings}
\usepackage{enumitem}
\usepackage[noend]{algpseudocode}
\usepackage{pdflscape}
\usepackage{afterpage}
\usepackage{array}
\usepackage{geometry}
\usepackage[square,numbers]{natbib}
\usepackage{verbatim}
\usepackage{epstopdf}
\usepackage{tikz}
\newcommand{\bx}{\mathbf{x}}

\newcommand{\be}{\mathbf{e}}

\newcommand{\R}{\mathbb{R}}

\newcommand{\bal}{\boldsymbol\alpha}

\input{macros.tex}
\input{tikz_setting.tex}

\DeclareMathOperator*{\minimize}{minimize}
\DeclareMathOperator*{\subjectto}{subject~to}

\DeclareMathOperator{\prox}{prox}
\DeclareMathOperator{\conv}{conv}
\DeclareMathOperator{\proj}{proj}

\newtheorem{definition}{Definition}
\newtheorem{example}{Example}
\newtheorem{theorem}{Theorem}
\newtheorem*{properties*}{Properties}

\usepackage{color}

\begin{document}

\title{A Primer on Coordinate Descent Algorithms}

\author{Hao-Jun Michael Shi\\
Department of Industrial Engineering and Management Sciences\\
Northwestern University\\
hjmshi@u.northwestern.edu
\and
 Shenyinying Tu\\
Department of Industrial Engineering and Management Sciences\\
Northwestern University\\
ShenyinyingTu2021@u.northwestern.edu
\and
Yangyang Xu\\
Department of Mathematics\\
University of Alabama\\
yangyang.xu@ua.edu
\and
Wotao Yin\\
Department of Mathematics\\
University of California, Los Angeles\\
wotaoyin@math.ucla.edu}


\maketitle
\setcounter{page}{2}
\tableofcontents


\begin{abstract}
\setcounter{page}{4}
This monograph presents a class of algorithms called coordinate descent algorithms for mathematicians, statisticians, and engineers outside the field of optimization. This particular class of algorithms has recently gained popularity due to their effectiveness in solving large-scale optimization problems in machine learning, compressed sensing, image processing, and computational statistics. Coordinate descent algorithms solve optimization problems by successively minimizing along each coordinate or coordinate hyperplane, which is ideal for parallelized and distributed computing. Avoiding detailed technicalities and proofs, this monograph gives relevant theory and examples for practitioners to effectively apply coordinate descent to modern problems in data science and engineering.\\

\noindent To keep the primer up-to-date, we intend to publish this monograph only after no additional topics need to be added and we foresee no further major advances in the area.
\end{abstract}

\setcounter{page}{5}
\section{Introduction} \label{introduction}

\subsection{Overview}

This monograph discusses a class of algorithms, called \textit{coordinate descent} (CD) algorithms, which is useful in solving large-scale optimization problems with smooth or non-smooth and convex or non-convex objective functions. Although these methods have existed since the early development of the discipline and the optimization community did not emphasize them until recently, 
various modern applications in machine learning, compressed sensing, and large-scale computational statistics have yielded new problems well suited for CD algorithms. These methods are generally applicable to a variety of problems involving large or high-dimensional data sets since they naturally break down complicated optimization problems into simpler subproblems, which are easily parallelized or distributed. For some structured problems, CD has been shown to perform faster than traditional algorithms, such as gradient descent. 
In addition, CD is generally applicable to non-convex problems and are easier to understand than splitting methods such as the Alternating Direction Method of Multipliers in this aspect. Also, few assumptions are needed to prove convergence to minima for convex problems and stationary points for non-convex problems. In fact, certain CD variants have also been shown to converge for non-convex functions with fairly loose properties.

CD algorithms follow the universal approach to algorithmic, numerical optimization: solving an optimization problem by solving a sequence of simpler subproblems. Each iterate is found by fixing most components of the variable vector $\vx$ at their current values and approximately minimizing the objective function with the remaining chosen components. In this monograph, we will explore a variety of interesting variants, extensions, and applications of CD connected to many different topics in optimization, statistics, and applied mathematics.

\subsection{Formulations}

We will consider the following general unconstrained minimization problem:
\begin{equation}\label{prob1}
\minimize_{\vx} f(\vx) = f(x_1, \ldots, x_n)
\end{equation}
where $\vx=(x_1,\ldots,x_n)\in\RR^n$ and the function $f:\RR^n \rightarrow \RR$ is continuous. Further assumptions may be made on the structure of $f$, such as convexity, Lipschitz continuity, differentiability, etc., while discussing theoretical guarantees for specific algorithms.

In addition, we will consider the following structured problem:
\begin{equation}\label{prob2}
\minimize_{\vx} F(\vx) = f(\vx) + \sum_{i = 1}^n r_i(x_i)
\end{equation}
where 
$f$ is differentiable, and $r_i$'s are extended-valued and possibly nondifferentiable functions.
Problems appearing in many recent applications such as compressed sensing, statistical variable selection, and model selection can be formulated in the form of (\ref{prob2}). Since we allow each $r_i$ to be extended-valued, it can model constraints on $x_i$ by including an indicator function, as discussed in Appendix \ref{sec: extended value}. The function $r_i$ can also include certain regularization terms to promote the structure of solutions, such as sparsity and low-rankness. We will further generalize the coordinate separable function $r_i$'s to block separable ones in Section \ref{section:block cd}.

\subsection{Framework of Coordinate Descent}

The basic coordinate descent (CD) framework for \eqref{prob1} and \eqref{prob2} is shown in Algorithm \ref{cd}. At each iteration, we choose one component $x_{i_k}$ and adjust it by a certain update scheme while holding all other components fixed.

\begin{algorithm}
\caption{Coordinate Descent}\label{cd}
\begin{algorithmic}[1]
\State Set $k = 0$ and initialize $\vx^0 \in \RR^n$;
\Repeat
        \State Choose index $i_k \in \{1,2,...,n\}$;
        \State Update $x_{i_k}$ to $x_{i_k}^k$ by a certain scheme depending on  $x^{k-1}$ and $f$ or $F$;
        \State Keep $x_j$ unchanged, i.e., $x_j^k = x_j^{k - 1}$, for all $j \neq i_k$;
        \State Let $k = k + 1$;
\Until{termination condition is satisfied};
\end{algorithmic}
\end{algorithm}

\begin{figure}
\center
\includegraphics[scale=.8]{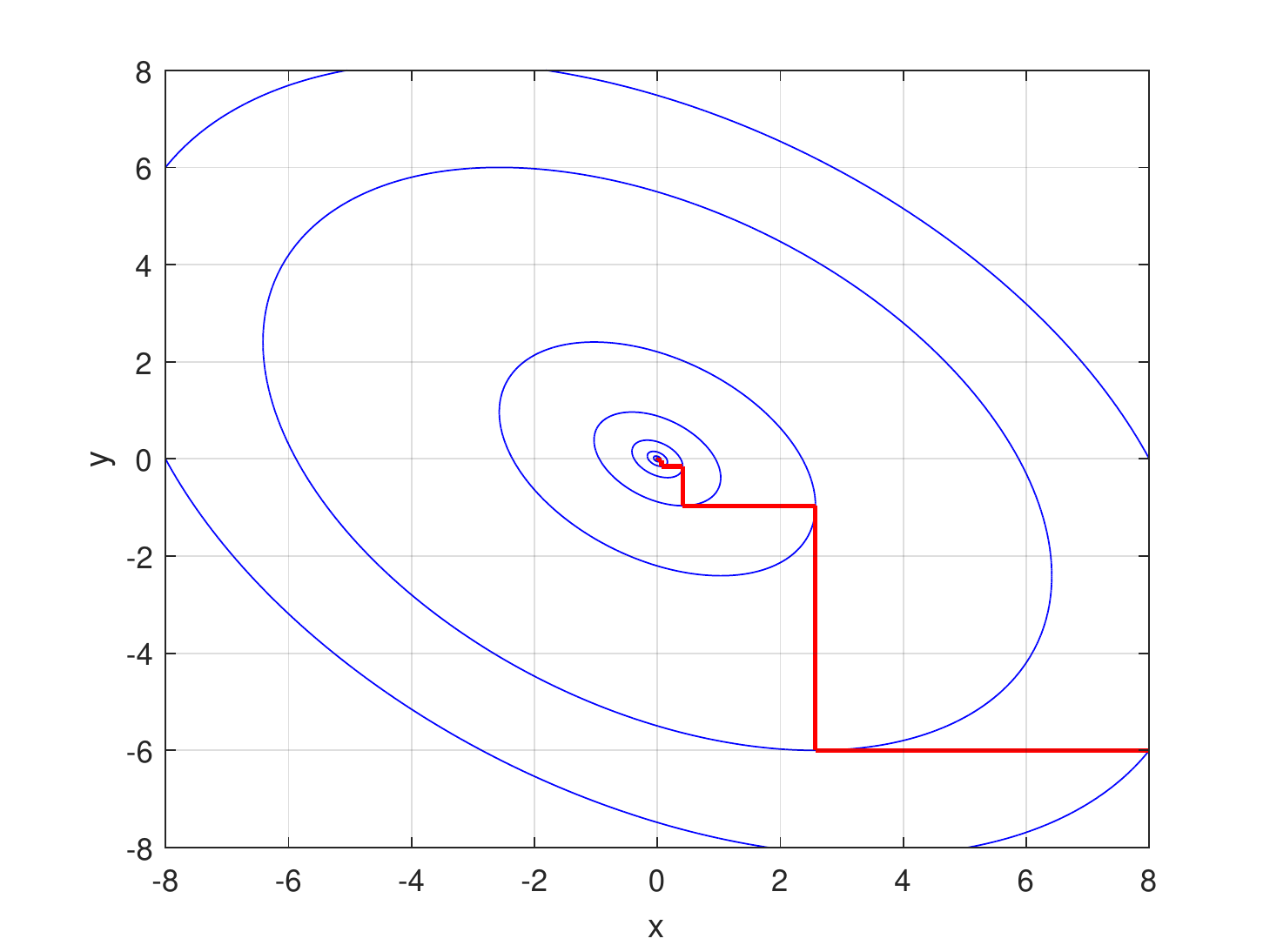}
\caption{CD applied on the quadratic function $f(x,y) = 7x^2 + 6xy + 8y^2$ with initial point $(8,-6)$. It minimizes $f$ alternatingly with respect to one of $x$ and $y$ while fixing the other. The blue curves correspond to different level curves of the function.}\label{2d}
\end{figure}

Intuitively, CD methods are easily visualized, particularly in the 2-dimensional case. Rather than moving all coordinates along a descent direction, CD changes a chosen coordinate at each iterate, moving as if it were on a grid, with each axis corresponding to each component. Figure \ref{2d} illustrates this process, where the coordinate minimization scheme is applied.

Within this framework, there are many different approaches for choosing an index and updating the selected coordinate. These various rules and updates affect the convergence properties of CD on different types of problems. They may exploit problem-specific structures, such as sparsity. They may also perform differently depending on the conditioning of the problem or how much each coordinate subproblem depends on one another.

The most natural approach to choosing an index is to select components cyclically, i.e. $i_0 = 1$, $i_1 = 2$, $i_2 = 3$, and so on. Alternatively, we can select a component at random at each iteration (not necessarily with equal probability). Lastly, we can choose components greedily, choosing the component corresponding to the greatest descent, strongest descent potential, or other scores, at the current iteration. The index rule may also satisfy an \textit{essentially cyclic} condition, in which every component is guaranteed to be updated at least once within every $N$ iterations. For example, we can perform a cyclic choice of components that are randomly shuffled after each cycle.

Regarding the update schemes, we can simply renew the selected component $x_{i_k}$ by minimizing the objective with respect to $x_{i_k}$ while fixing the remaining ones. Specifically, for problem \eqref{prob1}, we can perform the update:
\begin{equation}\label{cd-min}
x_{i_k}^k = \argmin_{x_{i_k}} f(x_1^{k - 1}, \ldots,x_{i_k - 1}^{k - 1}, x_{i_k}, x_{i_k + 1}^{k - 1}, \ldots, x_n^{k - 1}),
\end{equation}
and for problem \eqref{prob2}, one can have a similar update by replacing $f$ with $F$. Other update schemes can be performed if the problem has more structure. Suppose $f$ is differentiable in (\ref{prob1}), then one can apply coordinate-gradient descent along each chosen component, i.e.
\begin{equation}\label{classic1}
x_{i_k}^k = x_{i_k}^{k - 1} - \alpha_{i_k} \grad_{i_k} f(\vx^{k - 1})
\end{equation}
where $\alpha_{i_k}$ is a step size that can be set by line search or according to the property of $f$.  

The scheme in \eqref{classic1} can be easily extended to solving (\ref{prob2}) as 
\begin{equation}\label{classic2}
x_{i_k}^k = x_{i_k}^{k - 1} - \alpha_{i_k} (\grad_{i_k} f(\vx^{k - 1}) + \tilde{\grad} r_{i_k}(x^{k - 1}_{i_k}))
\end{equation}
where $\tilde{\grad} r_{i_k}(x_{i_k}^{k - 1})$ is a subgradient of $r_{i_k}$ at $x_{i_k}^{k - 1}$.
In addition, \textit{proximal} or \textit{prox-linear updates} can handle \eqref{prob2} when $r_i$'s are not differentiable. 
These updates minimize a surrogate function that dominates the original objective around the current iterate. The proximal update uses as the surrogate function the sum of the original function and a proximal term, and the prox-linear update employs a surrogate function to be the linearization of the differentiable part plus a proximal term and the nondifferentiable function. 
They both involve the \textit{proximal operator}, which for the function $\alpha f$ is defined as
$$\prox_{\alpha f} (\vy) = \argmin_{\vx} f(\vx) + \frac{1}{2 \alpha} \| \vx - \vy \|_2^2.$$
For many functions with favorable structures, the proximal operator is cheap to compute. We call these functions \emph{proximable functions}. Please refer to Appendix \ref{proximal operators} for more detail.

If we are minimizing the sum, or the average, of a vast number of functions, we can also apply \textit{stochastic updates}, which use sample gradients computed by selecting one or a few of these functions at each update instead of the exact gradient. We will discuss these alternative update schemes further in Section \ref{implementations}.

The step size $\alpha_{i_k}$ in \eqref{classic1} and \eqref{classic2} can also be chosen in many fashions. The choice of stepsize is important because a descent direction is not sufficient to guarantee descent. If $\alpha_{i_k}$ is too large, the function value may increase; if $\alpha_{i_k}$ is too small, the algorithm will converge at a relatively slow rate. To select the step size $\alpha_{i_k}$, we may perform an exact line search along the $i_k$th component, use traditional, backtracking line search methods to obtain sufficient descent, which may be more economical for certain separable objectives, or make predefined choices of $\alpha_{i_k}$ based on known properties of $f$.

Useful examples that shed light on the performance of CD on differently structured problems are given in Section \ref{applications}.
 We will also discuss related \textit{coordinate friendly} analysis and theory, as well as useful heuristics, applied to various problems in Section \ref{cf structures}.

\subsection{Other Surveys}

Coordinate descent algorithms have existed since the formation of the discipline. In response to the rising interest in large-scale optimization, a few articles have recently surveyed this class of algorithms. Wright \cite{wright2015coordinate} gives an in-depth review of coordinate descent algorithms, including convergence theory, which we highly recommend for optimizers and practitioners with a stronger background in optimization. Lange \cite{lange2014brief} provides a survey of optimization algorithms for statistics, including block coordinate descent algorithms.

We emphasize that this paper is specifically targeted towards engineers, scientists, and mathematicians outside of the optimization field, who may not have the requisite knowledge in optimization to understand research articles in this area. Because of this, we do not present the theory of coordinate descent algorithms in a formal manner but emphasize performance on real-world applications. In addition, we avoid discussing specific parallelized and distributed coordinate method implementations in detail since this remains an active area of research and conclusions cannot be drawn without discussing many implementation aspects. We instead give a list of possible approaches toward developing parallelized and distributed algorithms.  
We believe that the contents of this paper may serve as a guide and toolbox for practitioners to apply coordinate descent to more problems and applications in the future.

\subsection{Outline}

In Section \ref{implementations}, we give different classes of variants, including different update schemes, indexing schemes, as well as introduce the more generalized block coordinate formulation. In Section \ref{cf structures}, we give the relevant theory to analyze problems in practice for CD, and discuss useful heuristics to exploit \textit{coordinate friendly} structures. In Section \ref{applications}, we describe modern applications of CD in engineering and data science. In Section \ref{parallel}, we give resources for parallelizing CD for solving large-scale systems. Lastly, in Section \ref{conclusion}, we summarize our results from this monograph.

\subsection{Notation}

We introduce some notation before proceeding. 
Let $L$ be the gradient Lipschitz constant of the differentiable part $f$, i.e. for any $\vx, \vy$, it holds that
$$\| \grad f(\vx) - \grad f(\vy) \|_2 \leq L \|\vx - \vy\|_2. $$
Let $L_i$ denote the block-wise gradient Lipschitz constant, i.e. for any $\vx, \vy$, 
$$\| \grad_i f(\vx_1,\ldots,\vx_i,\ldots,\vx_s) - \grad_i f(\vx_1,\ldots,\vx_{i-1},\vy_i,\vx_{i+1},\ldots,\vx_s) \|_2 \leq L_i \| \vx_i - \vy_i \|_2,$$
where note that $L_i$ may depend on the value of $\vx_j$ for all $j\neq i$.
In addition we introduce the notation
$$f(\vx_{i_k}, \vx_{\neq i_k}) = f(\vx_1, ..., \vx_{i_k}, ..., \vx_s)$$
when the $i_k$th block is chosen.

\section{Algorithm Variants and Implementations}\label{implementations}

A wide range of implementation variants of CD have been developed for a large variety of applications. We discuss variants on the classic CD method introduced in Section \ref{introduction} and describe their strengths and weaknesses below.

\subsection{Block Coordinate Descent}\label{section:block cd}

Up until this point, we have only considered the method 
that updates one component of the variable $\vx$ at each iterate. We may generalize these coordinate updates to block coordinate updates. 
This method is particularly useful for applications with variables partitioned into blocks, such as non-negative matrix/tensor factorization (e.g., see \cite{cichocki2009nonnegative}), group LASSO \cite{yuan2006model}, and many distributed computing problems, where blocks of variables naturally appear.

Consider the following optimization problem:
\begin{equation}\label{prob3}
\minimize_{\vx} F(\vx)=f(\vx_1, ..., \vx_s) + \sum_{i = 1}^s r_i (\vx_i)
\end{equation}
where $\vx \in \RR^n$ is decomposed into $s$ block variables $\vx_1, ..., \vx_s$, 
$f$ is differentiable, and $r_i$ for $i = 1, ..., s$ are extended-valued and possibly nondifferentiable functions.
Note that if we consider the block formed by each component, we obtain formulation (\ref{prob2}). We may similarly adjust formulation (\ref{prob1}) for block variables.

From now on, we consider formulation (\ref{prob3}) since it is a generalization of formulation (\ref{prob2}). We can simply modify Algorithm \ref{cd} to fit this block structured problem, and the modified method is dubbed as Block Coordinate Descent (BCD), given in Algorithm \ref{bcd}.

\begin{algorithm}
\caption{Block Coordinate Descent}\label{bcd}
\begin{algorithmic}[1]
\State Set $k = 0$ and choose $\vx^0 \in \RR^n$;
\Repeat
        \State Choose index $i_k \in \{1,2,...,s\}$;
        \State Update $\vx_{i_k}$ to $\vx_{i_k}^k$ by a certain scheme depending on $\vx^{k-1}$ and $F$;
        \State Keep $\vx_j^k = \vx_j^{k - 1}$ for $j \neq i_k$;
        \State Let $k = k + 1$
\Until{termination condition is satisfied};
\end{algorithmic}
\end{algorithm}

Rather than updating a chosen coordinate at each iterate, block CD seeks to renew a chosen block of coordinates 
while other blocks are fixed. This method lends itself well for distributed or parallel computing since the update of a block coordinate is typically cheaper than that of all block variables. This will be discussed further in Section \ref{cf structures}. 

Block CD is also a generalization of the alternating minimization method that has been applied to a variety of problems, and also the expectation-maximization (EM) algorithm \cite{dempster1977maximum}, that performs essentially a 2-block CD. 

\subsection{Update Schemes}

As that done in \eqref{cd-min}, one can simply update $\vx_{i_k}$ by minimizing $F$ with respect to $\vx_{i_k}$ while fixing the remaining block variables. However, 
this update scheme can be hard since its corresponding subproblem may be difficult to solve exactly. In addition, BCD with this update scheme may not converge for some non-smooth and/or non-convex problems. This deficiency motivates the introduction of alternative update schemes 
that may give easier subproblems and ensure the convergence of the algorithm. 

All of the update schemes are summarized below:

\begin{enumerate}
\item \textit{(Block) Coordinate Minimization}:
$$\vx_{i_k}^k = \argmin_{\vx_{i_k}} f(\vx_{i_k}, \vx_{\neq i_k}^{k - 1}) + r_{i_k}(\vx_{i_k});
$$
\item \textit{(Block) Proximal Point Update}:
$$\vx_{i_k}^k = \argmin_{\vx_{i_k}} f(\vx_{i_k}, \vx_{\neq i_k}^{k - 1}) + \frac{1}{2\alpha_{i_k}^{k - 1}}\|\vx_{i_k} - \vx_{i_k}^{k - 1}\|_2^2 + r_{i_k}(\vx_{i_k});$$
\item \textit{(Block) Proximal Linear Update (Prox-Linear)}:
$$\vx_{i_k}^k = \argmin_{\vx_{i_k}} f(\vx^{k - 1}) +  \langle \grad_{i_k} f(\vx_{i_k}^{k - 1}, \vx_{\neq i_k}^{k - 1}), \vx_{i_k} - \vx_{i_k}^{k - 1}\rangle  + \frac{1}{2\alpha_{i_k}^{k - 1}}\|\vx_{i_k} - \vx_{i_k}^{k - 1}\|_2^2 + r_{i_k}(\vx_{i_k});$$
\end{enumerate}
where in the proximal point update, the step size $\alpha_{i_k}^{k-1}$ can be any bounded positive number, and in the prox-linear update, the step size $\alpha_{i_k}^{k-1}$ can be set to $1/L_{i_k}^{k - 1}$.


Since each update scheme solves a different subproblem, the updates may generate different sequences that converge to different solutions. The coordinate minimization, proximal point, and prox-linear updates may be interpreted as minimizing a surrogate function that upper bounds the original objective function when $\alpha_{i_k}$ is chosen appropriately, as noted in the BSUM algorithm \cite{hong2016unified}. It is important to understand the nuances of each scheme to apply the correct variant for a given application. We describe each update scheme in more detail below.

\subsubsection{Block Coordinate Minimization}

A natural way to update the selected block $\vx_{i_k}$ is to minimize the objective with respect to $\vx_{i_k}$ with all other blocks fixed, i.e., by the block coordinate minimization scheme:
\begin{equation}\label{classic}
\vx_{i_k}^k = \argmin_{\vx_{i_k}} f(\vx_{i_k}, \vx_{\neq i_k}^{k - 1}) + r_{i_k}(\vx_{i_k}).
\end{equation}

This classic scheme is most intuitive and was first introduced in \cite{hildreth1957quadratic} and further analyzed in \cite{beck2013convergence,d1959convex,grippo2000convergence,luo1992convergence,tseng2001convergence,warga1963minimizing}. 
BCD with this scheme is guaranteed to converge to a stationary point when the objective is convex, continuously differentiable, and strictly convex on each coordinate.
However, BCD may not converge for some nonconvex problems. Powell \cite{powell1973search} gives an example for which cyclic BCD with the update scheme in \eqref{classic} fails to converge to a stationary point. Although the block minimization scheme is most easily accessible and intuitive, alternative update schemes can have greater stability, better convergence properties, or easier subproblems. For Powell's example, BCD with the proximal point or prox-linear update schemes does converge.

\begin{figure}
\center
\includegraphics[scale = 0.5]{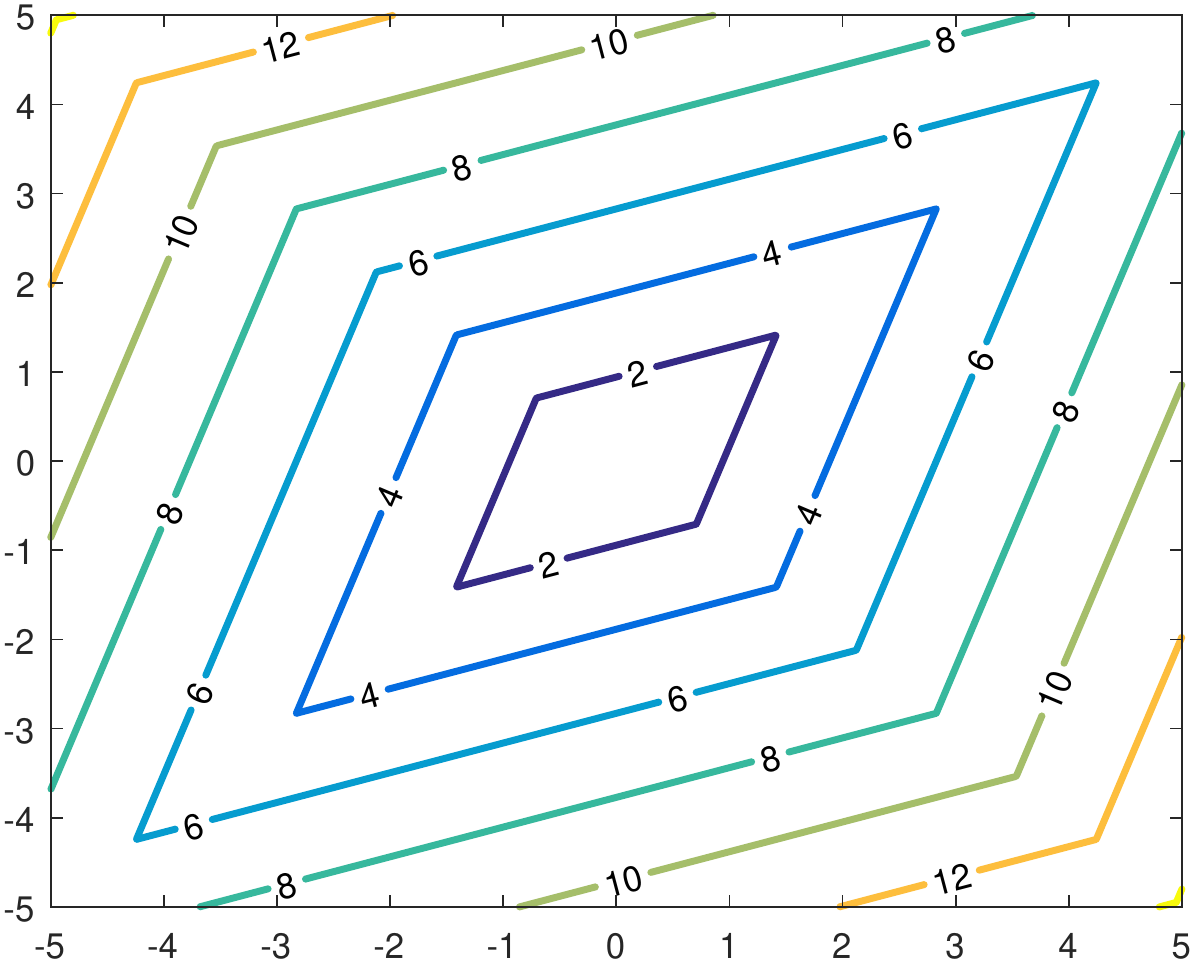}\qquad
\includegraphics[scale = 0.5]{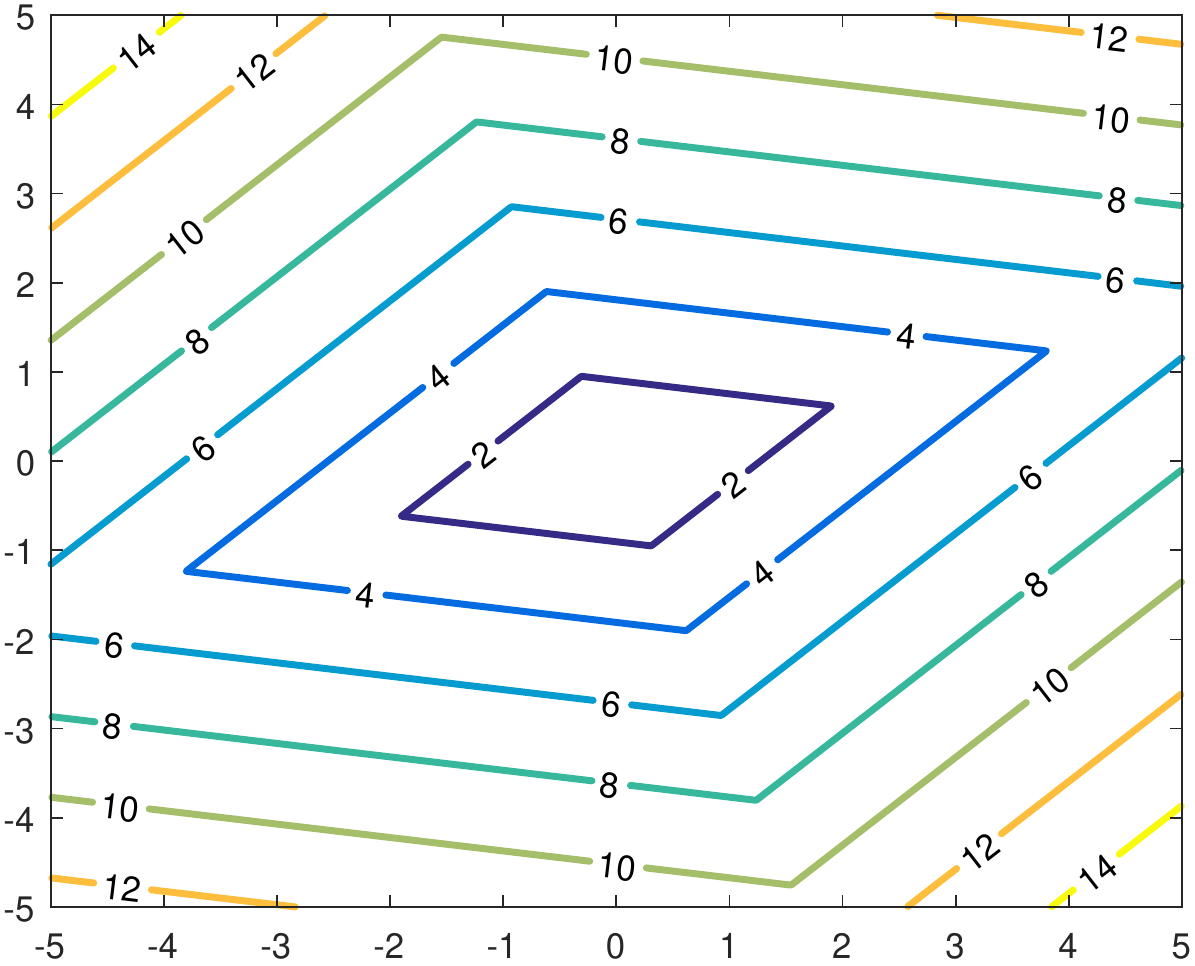}
\caption{The function $\ell(x,y)=|x|+2|y|$ is a separable, non-smooth, convex function. The left side is a $\pi/4$-radian rotation of $\ell$ for which CD gets stuck at a non-stationary point using any of the presented updates. The right side is a $\pi/10$-radian rotation of $\ell$ for which CD can correctly minimize.}\label{warga}
\end{figure}

Warga \cite{warga1963minimizing} provides a convex but non-smooth example, $f(x, y) = |x - y| - \min(x, y)$, for which BCD will get stuck at a non-stationary point. We illustrate this issue by rotating the simple function
$$\ell(x,y)=|x|+2|y|.$$
Figure \ref{warga} depicts the $\pi/4$ and $\pi/10$ radian rotations of $\ell$. If we consider the $\pi/4$ radian rotation of $\ell$ and start at a point of the form $(\beta, \beta)$ where $\beta \neq 0$, then coordinate minimization along $x$ or $y$ will not decrease the objective value since the point is already optimal along each coordinate direction. Since the optimal solution is at $(0, 0)$, the algorithm will get stuck at a non-optimal point. However, this problem does not occur for the $\pi/10$ radian rotation of $\ell$, since it can decrease along the $x$ direction.

More mathematically, the rotated functions correspond to setting $\varepsilon$ to $\pi/4$ and $\pi/10$, respectively, in the following function:
\begin{align*}
f_\varepsilon(x,y) = \ell\big(\cos(\varepsilon)x+\sin(\varepsilon)y,\cos(\varepsilon)y-\sin(\varepsilon)x\big).
\end{align*}
We can verify that $f_{\pi/4}(x,x)=\min_{\bar{y}} f_{\pi/4}(x,\bar{y})$ for any $x\in\RR$, that is, given any $x$, the minimizer $y$ of $f_{\pi/4}$ equals the value of $x$. (The same holds if, instead, we fix $y$ and minimize $f_{\pi/4}$ over $x$. The minimizer $x$ will equal the value of $y$.) Therefore, starting from any point $(x,y)$, a coordinate minimization over $x$ or $y$ will move to the point $(y,y)$ or $(x,x)$, respectively. If the point is not 0 (the origin), then any further coordinate update cannot reduce the value of $f_{\pi/4}$. Therefore, CD converges to a non-stationary point. 

Though our example $f_{\pi/4}$ gives a case where CD may converge to a non-stationary point if $f$ is convex, non-separable, and non-smooth, this may not always be the case. CD will converge for the example $f_{\pi/10}$ since it will not get stuck at any of the contour corners.

A mathematical explanation of this phenomenon is that the componentwise subgradients $p_x\in\partial_x f_{\pi/4}(x,y)$ and $p_y\in\partial_y f_{\pi/4}(x,y)$ do \emph{not} necessarily mean that the full vector formed by the concatenation of the component subgradients is a subgradient, i.e. $
[p_x;p_y]\in\partial f_{\pi/4}(x,y)$. Thus, in the case for $f_{\pi/4}$, a point in the form of $(\beta,\beta)\neq 0$ is not a stationary point because $0\not\in\partial f_{\pi/4}(\beta,\beta)$. A further explanation on this and subdifferential calculus is detailed in Appendix \ref{section: subdifferential calculus}. 


In general, BCD can fail to converge for an objective function that contains just one non-separable and non-smooth term, even if all the other terms are differentiable or separable.  In fact, this failure can occur with any BCD update presented in this paper. However, in this case, though there are no theoretical guarantees, practitioners may still try applying BCD and check optimality conditions for convergence. When the non-separable and non-smooth term has the form $f(Ax)$ and is convex, primal-dual coordinate update algorithms may also be applied since they will decouple $f$ from $A$, though this is not a focus in this monograph. Please refer to \cite{PesquetRepetti2014_class,peng2016coordinate,gao2016randomized} for more information on primal-dual coordinate update algorithms.

\subsubsection{Proximal Point Update}

The \textit{proximal point} update, or \textit{proximal} update, is defined by the following:
\begin{equation}\label{proximal}
\vx_{i_k}^k = \argmin_{\vx_{i_k}} f(\vx_{i_k}, \vx_{\neq i_k}^{k - 1}) + \frac{1}{2\alpha_{i_k}^{k - 1}}\|\vx_{i_k} - \vx_{i_k}^{k - 1}\|_2^2 + r_{i_k}(\vx_{i_k}),
\end{equation}
where $\alpha_{i_k}^{k-1}$ serves as a step size and can be any bounded positive number.

The proximal update with BCD was introduced in \cite{auslender1992asymptotic} for convex problems and further analyzed in \cite{grippo2000convergence,razaviyayn2013unified,xu2013block} for possibly nonconvex problems. This update adds a quadratic proximal term to the classic block minimization update described in \eqref{classic}, and thus the function of each subproblem dominates the original objective around the current iterate. This modification gives the proximal update scheme better convergence properties and increased stability, particularly for non-smooth problems.

\subsubsection{Prox-Linear Update}

The \textit{proximal linear} update, or \textit{prox-linear} update, is defined by:
\begin{equation}\label{prox-linear}
\vx_{i_k}^k = \argmin_{\vx_{i_k}} f(\vx^{k - 1}) + \langle \grad_{i_k} f(\vx_{i_k}^{k - 1}, \vx_{\neq i_k}^{k - 1}), \vx_{i_k} - \vx_{i_k}^{k - 1}\rangle  + \frac{1}{2\alpha_{i_k}^{k - 1}}\|\vx_{i_k} - \vx_{i_k}^{k - 1}\|_2^2 + r_{i_k}(\vx_{i_k}),
\end{equation}
where the step size $\alpha_{i_k}^{k-1}$ can be set as the reciprocal of the Lipschitz constant of $\nabla_{i_k} f(\vx_{i_k},\vx_{\neq i_k}^{k-1})$ with respect to $\vx_{i_k}$. The difference between the prox-linear update from the proximal update is that the former further linearizes the smooth part $f$ to make the subproblem easier.

The prox-linear update scheme was introduced in \cite{tseng2009coordinate}, which proposed a more general framework of block coordinate gradient descent (BCGD) methods. It was later adopted and also popularized by Nesterov \cite{nesterov2012efficiency} in the randomized coordinate descent method. BCD with this update scheme has been analyzed and also applied to both convex and nonconvex problems such as in \cite{beck2013convergence, hong2013iteration, nesterov2012efficiency, tseng2009coordinate, xu2014globally, yun2011block, yun2011coordinate, zhou2015global}. 
In essence, this scheme minimizes a surrogate function that 
dominates the original objective around the current iterate $\vx_{i_k}^{k - 1}$. 
Note that when the regularization term $r_{i_k}(\vx_{i_k})$ vanishes, this update reduces to:
\begin{equation}\label{prox-linear without r}
\vx_{i_k}^k = \argmin_{\vx_{i_k}} f(\vx^{k - 1}) + \langle \grad_{i_k} f(\vx_{i_k}^{k - 1}, \vx_{\neq i_k}^{k - 1}), \vx_{i_k} - \vx_{i_k}^{k - 1}\rangle  + \frac{1}{2\alpha_{i_k}^{k - 1}}\|\vx_{i_k} - \vx_{i_k}^{k - 1}\|_2^2,
\end{equation}
or equivalently the block gradient descent algorithm:
$$\vx_{i_k}^k= \vx_{i_k}^{k - 1}-\alpha_{i_k}^{k - 1}\grad_{i_k} f(\vx_{i_k}^{k - 1}, \vx_{\neq i_k}^{k - 1}).$$

When BCD with the proximal or prox-linear update scheme is applied to non-convex objectives, it often gives solutions of lower objective values compared with the block coordinate minimization, since small regions containing certain local minima may be avoided by its local proximal or prox-linear approximation. Note that the prox-linear update may take more iterations to reach the same accuracy than the other two schemes. However, it is easier to perform the update and thus may take less total time as demonstrated in \cite{shi2014sparse, xu2015alternating,xu2013block}.

\subsubsection{Extrapolation}\label{section:extrapolation}

Though the prox-linear update is easily computed and gives a better solution overall, coordinate minimization and proximal updates tend to make larger objective decreases per iteration. This observation motivates the use of extrapolation to accelerate the convergence of the prox-linear update scheme.

Extrapolation uses the information at an extrapolated point in place of the current point to update the next iterate \cite{beck2009fast, lee2013efficient, nesterov2013gradient, shalev2016accelerated, xu2013block}. In particular, rather than using the partial gradient at $\vx_{i_k}^{k - 1}$ for the next update, we instead consider an extrapolated point
\begin{equation}\label{extrapolation}
\hat{\vx}_{i_k}^{k - 1} = \vx_{i_k}^{k - 1} + \omega_{i_k}^{k - 1} (\vx_{i_k}^{k - 1} - \vx_{i_k}^{k - 2}),
\end{equation}
where $\omega_{i_k}^{k - 1} \geq 0$ is an extrapolation weight. This extrapolated point is then used to compute our next update, and it gives the update:
$$\vx_{i_k}^k = \argmin_{\vx_{i_k}} f(\vx^{k - 1}) + \langle \grad_{i_k} f(\hat{\vx}_{i_k}^{k - 1}, \vx_{\neq i_k}^{k - 1}), \vx_{i_k} - \hat{\vx}_{i_k}^{k - 1}\rangle  + \frac{1}{2\alpha_{i_k}^{k - 1}}\|\vx_{i_k} - \hat{\vx}_{i_k}^{k - 1}\|_2^2 + r_{i_k}(\vx_{i_k}).$$
Note that if $\omega_{i_k}^{k - 1}=0$, the above update reduces to that in \eqref{prox-linear}. Appropriate positive weight can significantly accelerate the convergence of BCD with prox-linear update as demonstrated in \cite{nesterov2012efficiency, xu2015alternating} while the per-iteration complexity remains almost the same.

\subsubsection{Stochastic Gradients}

Often in machine learning and other large data applications, we encounter problems with datasets consisting of tens of millions of datapoints and millions of features. Due to the size and scope of these problems, sometimes computing a coordinate gradient may still be extremely expensive. To remedy this, we can introduce \textit{stochastic gradients}.

Consider the following stochastic program, called the regularized \textit{expected risk minimization} problem,
\begin{equation}\label{stoch-prob}
\minimize_{\vx} ~\mathbb{E}_{\xi} f_\xi(\vx) + \sum_{i = 1}^s r_i(\vx_i),
\end{equation}
where $\xi$ is a random variable, $f(\vx)=\mathbb{E}_{\xi} f_\xi(\vx)$ is differentiable, and $r_i$'s are certain regularization terms. 

The function $f$ may represent some loss due to inaccurate predictions from some classifier or prediction function. To minimize loss for any set of predicted and true parameters, we take the expectation of the loss with respect to some probability distribution modeled by the random variable $\xi$. This expectation takes the form of an integral or summation that weights losses for all possible predictions and true parameters.

An interesting case of \eqref{stoch-prob} is when $\xi$ follows a uniform distribution over $1,\ldots,m$, representing a noninformative prior distribution. In this case, if the potential outcomes are discrete, then the stochastic program in \eqref{stoch-prob} reduces to the \textit{empirical risk minimization} problem
\begin{equation}\label{sum of functions}
\minimize_{\vx} ~ \frac{1}{m} \sum_{i = 1}^m f_i(\vx) + \sum_{i = 1}^s r_i(\vx_i)
\end{equation}
where each $f_i$ represents the loss incurred with respect to one sample from the dataset. Therefore, since the problems we are considering often rely on millions of training points, $m$ is very large. The empirical risk minimization problem is also often used in place of the expected risk minimization problem when the probability distribution of $\xi$ is unknown.

The stochastic gradient method, also called the stochastic approximation method, (e.g., see \cite{nemirovski2009robust}) is useful for minimizing objectives in the form of \eqref{stoch-prob} or \eqref{sum of functions} with large $m$, in which computing the exact gradient or objective becomes overly expensive. To compute the full gradient for \eqref{sum of functions}, one would have to compute
$$\grad f(\vx) = \frac{1}{n} \sum_{i = 1}^n \grad f_i(\vx)$$
by processing every training example in the dataset. Since this is often infeasible, we can instead sample either one or a small batch of loss functions $f_i$, called a \textit{mini-batch}, to compute a subsampled gradient to use in place of the full gradient, i.e. we use
$$\tilde{\vg}_{i_k} = \frac{1}{|S_k|} \sum_{l \in S_k} \grad_{\vx_{i_k}} f_{k_l}(\vx)$$
where $S_k \subset \{1, ..., m\}$ is a mini-batch and $|S_k|$ is the number of sample functions selected from the loss functions $f_i$'s.

More generally, for solving \eqref{stoch-prob} or \eqref{sum of functions} by prox-linear BCD, we may replace the true coordinate gradient with a stochastic approximation, i.e. if the $i_k$th block is selected,
\begin{equation}\label{stoch-bcd}\vx_{i_k}^k = \argmin_{\vx_{i_k}} f(\vx^{k - 1}) + \langle \tilde{\vg}_{i_k}^{k - 1}, \vx_{i_k} - \vx_{i_k}^{k - 1} \rangle + \frac{1}{2\alpha_{i_k}^{k - 1}} \| \vx_{i_k} - \vx_{i_k}^{k - 1}\|^2 + r_{i_k}(\vx_{i_k}),
\end{equation}
where $\tilde{\vg}_{i_k}^{k - 1}$ is a stochastic approximation of $\nabla_{i_k}f(\vx^{k-1})$, $\tilde{\vg}_{i_k}^{k - 1}$ is a subsampled gradient, etc. 

Though the stochastic prox-linear update may be less accurate, it works well when there is a limited amount of memory available, or when a solution is needed quickly, as discussed in \cite{dang2013stochastic,xu2015block}. 

\subsubsection{Variance Reduction Techniques}

Alternatively, we can also consider \textit{stochastic variance-reduced gradients}, which use a combination of stale gradients with new gradients to reduce the variance in the chosen stochastic gradients. These variance-reduced stochastic gradient algorithms gain a significant speedup in convergence; whereas stochastic gradients only have sublinear convergence rate guarantees, variance-reduced stochastic gradients can have linear convergence rates similar to traditional gradient descent methods on problems with strongly convex objectives.

Consider the problem given above in (\ref{sum of functions}). Let $\phi_i^k$ denote the past stored point used at the prior gradient evaluation for function $f_i$ and $\grad f_i (\phi_i^k)$ denote the stored gradient. We list some common stochastic variance-reduced gradients below:
\begin{itemize}
\item SAG \cite{schmidt2013minimizing}: If the $j$th indexed function is chosen at iterate $k$,
$$\tilde{\vg}^{k - 1} = \frac{\grad f_j(\vx^{k - 1}) - \grad f_j(\phi_j^{k - 1})}{m} + \frac{1}{m} \sum_{l = 1}^{m} \grad f_l (\phi_l^{k - 1}).$$
The current iterate $\vx^{k - 1}$ is then taken as $\phi_j^k$ and $\grad f_j (\phi_j^{k})$ is explicitly stored in a table of gradients.
\item SAGA \cite{defazio2014saga}: If the $j$th indexed function is chosen at iterate $k$,
$$\tilde{\vg}^{k - 1} = \grad f_j(\vx^{k - 1}) - \grad f_j(\phi_j^{k - 1}) + \frac{1}{m} \sum_{l = 1}^{m} \grad f_l (\phi_l^{k - 1}).$$
The current iterate $\vx^{k - 1}$ is then taken as $\phi_j^k$ and $\grad f_j (\phi_j^{k})$ is explicitly stored in a table of gradients.
\item SVRG \cite{johnson2013accelerating}: If the $j$th indexed function is chosen at iterate $k$,
$$\tilde{\vg}^{k - 1} = \grad f_j(\vx^{k - 1}) - \grad f_j(\tilde{\vx}) + \frac{1}{m} \sum_{l = 1}^{m} \grad f_l (\tilde{\vx}),$$
where $\tilde{\vx}$ is not updated every step but is updated after a fixed number of iterations. If enough memory is available, individual gradients $\grad f_l (\tilde{\vx})$'s as well as their average are all stored; otherwise, one can store only the average and evaluate $\nabla f_j(\tilde \vx)$ at each iteration (in addition to the usual work to evaluate $\nabla f_j(\vx^{k-1})$). 
\end{itemize}

Another approach, the stochastic dual coordinate ascent algorithm (SDCA) \cite{shalev2013stochastic}, applies randomized dual coordinate ascent to the dual formulation of the problem and gives similar variance reduction properties.

In general, though these methods give better convergence properties, they require more memory to store stale gradients or more computation to evaluate exact gradient. However, they perform better than traditional stochastic gradient methods, and work well when calculating the exact gradient is expensive.

Note that the primary difference between SVRG and SAG is that SVRG makes 2-3x more gradient evaluations if it does not store a table of gradients, whereas SAG uses less gradient evaluations but requires more memory overhead to store gradients. SAGA may be interpreted as the midpoint between SVRG and SAG. The usage of SAG, SAGA, and SVRG is, therefore, problem-dependent.

The variance reduction technique can also be incorporated into the prox-linear BCD. For example, one can use any $\tilde{\vg}^{k-1}$ of the above three ones in the update \eqref{stoch-bcd} to accelerate its convergence.

%

\subsubsection{Summative Proximable Functions}

We apply the prox-linear update \eqref{prox-linear} when the function $r_{i_k}$ is proximable, that is, when its proximal operator can be evaluated at a low cost.
Functions such as $\ell_1$-norm, $\ell_2$-norm, and $\ell_\infty$-norm, as well as the indicator functions of box constraints, one or two linear constraints, and the standard simplex, are proximable functions.
The list can be quite long. Nonetheless, it is not difficult to see that, even if two functions $f$ and $g$ are both proximable, $f+g$ may not be proximable.
Therefore, the update \eqref{prox-linear} can still be expensive to compute if $r_{i_k}$ is the sum of two or more proximable functions.

The \emph{summative proximable function} is the sum of proximable functions $f$ and $g$ that satisfy $$\prox_{f+g}=\prox_{g}\circ\prox_{f}.$$ Because their proximal operator can be obtained by applying the proximal operator of $f$ and then that of $g$ in a sequential fashion, it is also proximable. Some common examples of summative proximable functions are:
\begin{itemize}
  \item $f(\vx)+g(\vx):=f(\vx)+\beta\|\vx\|_2$, where $\beta\ge 0$ and $f(\vx)$ is a homogeneous function of order 1 (i.e., $f(\alpha \vx) = \alpha f(\vx)$ for $\alpha\ge 0$). Examples of $f(\vx)$ include $\alpha\|\vx\|_1$, $\alpha\|\vx\|_\infty$, $\iota_{\ge 0}(\vx)$, or $\iota_{\le 0}(\vx)$ and $\alpha,\beta>0$.
  \item $f(\vx)+g(\vx):=\beta\mathrm{TV}(\vx)+g(\vx)$, where $\mathrm{TV}(\vx):=\sum_{i=1}^{n-1}|x_{i+1}-x_i|$ is (discrete) total variation, and $g(\vx)$ is  a function with the following property: for any $\vx\in\RR^n$ and coordinates $i\in[n]$ and $j=i+1$,
      \begin{align*}
        x_i>x_j&\quad\Rightarrow\quad\big(\prox_{g}(\vx)\big)_i \ge \big(\prox_{g}(\vx)\big)_j \\
        x_i<x_j&\quad\Rightarrow\quad\big(\prox_{g}(\vx)\big)_i \le \big(\prox_{g}(\vx)\big)_j \\
        x_i=x_j&\quad\Rightarrow\quad\big(\prox_{g}(\vx)\big)_i  =  \big(\prox_{g}(\vx)\big)_j.
      \end{align*}
Examples of such $g(\vx)$ include $\alpha\|\vx\|_1$, $\alpha\|\vx\|_2$, $\alpha\|\vx\|_\infty$, $\iota_{\ge 0}(\vx)$, or, more generally, $\iota_{[\ell,u]}(\vx)$ for any $\ell,u\in\RR$.
  \item \cite[Prop. 3.6]{CombettesPesquet2007_proximal} scalar function $f_i(\rho)+g_i(\rho):=\alpha|\rho|+g_i(\rho)$, where $\rho\in\RR$ and $g_i$ is convex and $g_i'(0)=0$. An example is the \emph{elastic net} regularizer \cite{ZouHastie2005_regularization}: $f(\vx)+g(\vx):=\alpha\|\vx\|_1+\frac{1}{2}\|\vx\|_2^2$.
\end{itemize}
The key to these results is an inclusion property: For any $\vx\in\RR^n$, let $\vy:=\prox_f(\vx)$ and $\vz:=\prox_g(\vy)$, whose minimization conditions are
\begin{align*}
  0  & \in \partial f(\vy)+(\vy-\vx), \\
  0  & \in \partial g(\vz)+(\vz-\vy),
\end{align*}
respectively, and adding them yields
\begin{align*}
  0  & \in \partial f(\vy)+\partial g(\vz)+(\vz-\vx).
\end{align*}
If the property of $f$ and $g$ gives the inclusion property $\partial f(\vy)\subseteq \partial f(\vz)$, then we arrive at the minimization condition of $\vz=\prox_{f+g}(\vx)$:
\begin{align*}
  0  & \in \partial f(\vz)+\partial g(\vz)+(\vz-\vx).
\end{align*}

Because the first two classes of summative proximable functions are not seen elsewhere to the best of our knowledge, a proof is included in Appendix \ref{proofs}.

\subsection{Choosing Update Index $i_k$}

In this section, we elaborate on various implementation approaches in choosing the coordinate or block $i_k \in \{1, ..., s\}$. Since different paths taken in coordinate descent may lead to different minima and different schemes perform differently for both convex and non-convex problems, the choice of the update index $i_k$ for each iterate is crucial for good performance for BCD. Often, it is easy to switch index orders. However, the choice of index affects convergence, possibly resulting in faster convergence or divergence. 
We describe the index rules more in detail below.

\subsubsection{Cyclic Variants}

The most natural, deterministic approach for choosing an index is to choose indices in a cyclic fashion, i.e. $i_0 = 1$ and
$$i_{k + 1} = (k\mod{s}) + 1,\ k \in \mathbb{N}.$$

We may also adapt this method and instead cycle through a permutation of $\{1, ..., s\}$, called a \textit{shuffled cyclic} method. In practice, one may reshuffle the order of the indices after each cycle, or cycle through all coordinates, which may have stronger convergence properties for some applications.

Another approach is to satisfy an \textit{essentially cyclic} condition, in which for every consecutive $N \geq s$ iterations, each component is modified at least once. 
More rigorously, we require
$$\bigcup_{j = 0}^N \{i_{k - j}\} = \{1, 2, ..., s\}$$
for all $k \geq N$.

Cyclic variants are most intuitive and easily implemented. BCD with the deterministic cyclic rule may give poorer performance than that with shuffled cyclic one, as demonstrated in \cite{xu2014globally} for solving non-negative matrix factorization. 
Convergence results of cyclic BCD are given in \cite{beck2013convergence, bonettini2011inexact, d1959convex, grippo2000convergence, hong2013iteration,
luo1992convergence, razaviyayn2013unified, saha2013nonasymptotic,
tseng2009coordinate, xu2014globally, zadeh1970note}. 

\subsubsection{Randomized Variants}

In randomized BCD algorithms, the update component or block $i_k$ is chosen randomly at each iteration. The simplest, most commonly used randomized variant is to simply select $i_k$ with equal probability, or \textit{sample uniformly}, independent of all choices made in prior iterations.

Other typical randomized variants include sampling without replacement, and considering different, non-uniform probability distributions. We list common sampling methods below:
\begin{enumerate}
\item \textit{Uniform sampling} \cite{fercoq2015accelerated, nesterov2012efficiency,richtarik2014iteration,shalev2011stochastic,shalev2013stochastic}: Each block coordinate $j \in \{1, ... s\}$ is chosen with equal probability as we described above, i.e.
$$P(i_k = j) = \frac{1}{s},\, j=1,\ldots,s.$$

\item \textit{Importance sampling} \cite{leventhal2010randomized,nesterov2012efficiency, richtarik2014iteration, zhang2004solving, zhang2014stochastic}: We proportionally weight each block according to its block-wise Lipschitz gradient constant $L_j > 0$ for all $j \in \{1, ..., s\}$.
More rigorously, given some $\alpha \ge 0$, 
Nesterov \cite{nesterov2012efficiency} proposes the following distribution:
$$P(i_k = j) = p_{\alpha}(j) = \dfrac{L_j^{\alpha}}{\sum_{i = 1}^s L_i^{\alpha}},\ j=1,\ldots, s.$$

This scheme generalizes uniform sampling -- when $\alpha = 0$, we obtain uniform sampling. Importance sampling has been further studied in \cite{Allen-ZhuQuRichtarikYuan2015_EvenFaster,CsibaQuRichtarik2015_StochasticDual}.

\item \textit{Arbitrary sampling} \cite{patrascu2015efficient,qu2014coordinate1, qu2014coordinate2,richtarik2014iteration, richtarik2015optimal}: We pick and update a block $j \in \{1, 2, ..., s\}$ arbitrarily, following some assigned probability distribution $(p_1, ..., p_s)$, i.e.
$$P(i_k = j) = p_j,\ j=1,\ldots,s,$$
where $0 \leq p_i \leq 1$ for all $i$ and $\sum_{i = 1}^s p_i = 1$. This sampling scheme generalizes both uniform and importance sampling.
\end{enumerate}

Randomized BCD is well suited for cases in which memory or available data is limited since the computation of a partial derivative is often much cheaper and less memory demanding than computing an entire gradient \cite{nesterov2012efficiency}. Recent work also suggests that randomization improves the convergence rate of BCD in expectation, such as for minimizing generic smooth and simple nonsmooth block-separable convex functions \cite{nesterov2012efficiency, qu2014coordinate1, richtarik2014iteration,tao2012stochastic}. 
However, randomized BCD variants have greater per-iteration complexities than cyclic BCD variants since these algorithms have to sample from probability distributions each iteration. 
Since randomized variants are non-deterministic, results in practice may vary. During the running of a randomized BCD, cache misses are more likely, requiring extra time to move data from slower to faster memory in the memory hierarchy.

\subsubsection{Greedy Variants}\label{greedy}

The last widely used approach for selecting indices are greedy methods. These methods choose $i_k$ ``greedily", or choose the index such that the objective function is minimized most, or the gradient or subgradient has the largest size, in that direction. The simplest variant for smooth functions, called the \textit{Gauss-Southwell\footnote{The greedy selection rule dates back to Gauss and was popularized by Southwell \cite{Southwell1940_relaxation} for linear systems.} selection rule (GS)}, involves choosing $i_k$ such that the gradient of the chosen block is greatest, or mathematically,
\begin{equation}\label{GS}
i_k = \argmax_{1 \leq j \leq s} \| \grad_j f(\vx^{k - 1}) \|
\end{equation}
for formulation (\ref{prob1}).  
If each block consists of an individual component, then the norm reduces to the absolute value function. GS can be analyzed in the general CD framework \cite{tseng1990dual,luo1992convergence} for convex optimization. In the refined analysis \cite{nutini2015coordinate}, it is shown that, except in extreme cases, GS converges faster than choosing random coordinates in terms of the number of iterations, though its per-iteration complexity is higher. 

Alternatively, we could choose the component or block that gives maximum improvement, called the \textit{Maximum Block Improvement (MBI) rule} \cite{chen2012maximum,li2015convergence}:
\begin{equation}\label{MBI}
i_k = \argmin_j f(\vx_j, \vx_{\neq j}^{k - 1})
\end{equation}

Motivated by the performance of Lipschitz sampling and Gauss-Southwell's rule, Nutini, et. al. \cite{nutini2015coordinate} proposed the \textit{Gauss-Southwell-Lipschitz (GSL)} rule:
\begin{equation}\label{GSL}
i_k = \argmax_{j} \frac{\| \grad_j f(\vx^{k - 1}) \|}{\sqrt{L_j}}.
\end{equation}
The GSL rule accounts for varied coordinatewise Lipschitz constants. When the gradients of two coordinates have similar sizes, updating the coordinate that has the smaller $L_i$ will likely lead to greater reduction in the objective value and is thus preferred in the selection.

For non-smooth problems, we list some commonly used GS rules for formulation (\ref{prob3}). We let $L \in \RR$ denote the gradient Lipschitz constant and $L_j \in \RR$ denote the $j$th coordinate gradient Lipschitz constant.

\begin{enumerate}
\item \textit{Gauss-Southwell-s rule (GS-s)}: At each iteration, the coordinate $i_k \in \{1, ..., s\}$ is chosen by
\begin{equation}\label{gs-s}
i_k = \argmax_{j} \left\{ \min_{\tilde{\grad} r_j \in \partial r_j} \| \grad_j f(\vx^{k - 1}) + \tilde{\grad} r_j(\vx_j^{k - 1}) \| \right\}.
\end{equation}
This rule chooses the coordinate with the greatest negative partial derivative, similar to (\ref{GS}). It is popular in $\ell_1$ minimization \cite{li2009coordinate,shevade2003simple,WuLange2008_CoordinateDescent}.


\item \textit{Gauss-Southwell-r rule (GS-r)}: At each iteration, the coordinate $i_k \in \{1, ..., s\}$ is chosen by
\begin{equation}\label{gs-r}
i_k = \argmax_{j} \big\| \vx_j^{k - 1} - \prox_{\frac{1}{L} r_j} [\vx_j^{k - 1} - \frac{1}{L} \grad_j f(\vx^{k-1}) ] \big\|
\end{equation}
which is equivalent to
$$i_k = \argmax_{j} \big\| \argmin_{\vd} \big(f(\vx^{k - 1}) + \grad_j f(\vx^{k - 1})^T \vd + \frac{L}{2} \|\vd\|_2^2 + r_j(\vx_j^{k - 1} + \vd) - r_j(\vx_j^{k - 1})\big) \big\|,$$
where $\vd$ has the same size as the block gradient direction.
This rule chooses the block that maximizes the distance from the current iterate to the block's following iterate from a proximal gradient update. It has been used in \cite{tseng2009coordinate,dhillon2011nearest,peng2013parallel}.
If the coordinate-wise gradient Lipschitz constant $L_j$ is known, then we can replace $L$ by $L_j$ in the GS-r update to obtain the \textit{Gauss-Southwell-Lipschitz-r rule (GSL-r)}.

\item \textit{Gauss-Southwell-q rule (GS-q)}: At each iteration, the coordinate $i_k \in \{1, ..., s\}$ is chosen by
\begin{equation}\label{gs-q}
i_k = \argmin_j \big( \min_{\vd} f(\vx^{k - 1}) + \grad_j f(\vx^{k - 1})^T \vd + \frac{L}{2} \|\vd\|_2^2 + r_j(\vx_j^{k - 1} + \vd) - r_j(\vx_j^{k - 1})  \big).
\end{equation}
This rule can be interpreted as the maximum coordinate-wise descent and has been used in \cite{tseng2009coordinate}. If the coordinate-wise gradient Lipschitz constant $L_j$ is known, then we can replace $L$ by $L_j$ in the GS-q update to obtain the \textit{Gauss-Southwell-Lipschitz-q rule (GSL-q)}.
\end{enumerate}

Note that these greedy methods usually require evaluating the whole gradient vector or some other greedy score, and searching for the best index. These greedy scores may be stored and updated using a max-heap, which is further detailed in \cite{meshi2012convergence}. To practically implement greedy methods, some terms of these greedy scores may be cached and maintained at each iteration to make these approaches computationally worthwhile. Section \ref{cf structures} provides detailed further explanations.

Greedy coordinate selections are very efficient for sparse optimization since most zero components in the solution are never selected and thus remain zero throughout the iterations \cite{li2009coordinate,peng2013parallel}. The problem dimension effectively reduces the number of variables that are ever updated, which is relatively small. Consequently, the greedy CD iteration converges in very few iterations. The saved iterations over-weigh the extra cost of ranking the coordinates.


Other simple greedy methods involve perturbing coordinates or blocks $\vx_j$ for all $j \in \{1, ..., s\}$ by some small step size $\beta > 0$, then evaluating the difference of objective values at those perturbed points with the current point to determine the coordinate or block of steepest descent, i.e.,
$$i_k = \argmax_j |f(\vx_j + \beta \be_j) - f(\vx_j)|,$$
where $e_j$ is the standard basis vector consisting of 1 at the $j$th entry and 0's elsewhere.

\subsubsection{Comparison of Index Rules}

We summarize the strengths and weaknesses of these common index rules below in Table \ref{index rules}. We also point readers to Section \ref{applications} 
for a comparison of various index rules applied to several examples. 
Note that no matter which index rule is chosen, CD can stagnate at non-critical points for objectives with terms that are both non-separable and non-smooth.

\begin{landscape}
\begin{table}
\footnotesize
\center
\begin{tabular}{| c | p{1.5in} | p{2in} | p{2in} |}
\hline
Method & Description & Strengths & Weaknesses \\ \hline
Cyclic & $i_k$ is chosen in a cyclically: $$i_{k} = (k \mod{n}) + 1$$
or by a randomly shuffled cycle
 &
\begin{itemize}[nosep, leftmargin =*]
\item Fast, easy implementation
\item Lowest per-iteration complexity
\item Performs well for problems with low coupling between blocks
\end{itemize} &
\begin{itemize}[nosep, leftmargin =*]
\item When coordinates are highly coupled, performance may be worse than randomized and greedy methods and theoretical worst-case complexity is worse
\item When exact coordinate minimization is used on non-convex problems, the points may cycle and fail to converge
\end{itemize} \\ \hline
Randomized & $i_k$ is chosen randomly following some probability distribution given by the vector $(p_1, ..., p_s)$: $$P(i_k = j) = p_j$$ &
\begin{itemize}[nosep, leftmargin =*]
\item Easier to analyze; often reduces to standard case after taking expectations
\item Empirically avoids local solutions more for some non-convex problems
\item Well suited for parallel computing
\end{itemize} &
\begin{itemize}[nosep, leftmargin =*]
\item Non-deterministic
\item Random data moves are slower and result in cache misses
\item Slightly higher per-iteration complexity than cyclic CD due to pseudo-random number generation
\item 2-3x slower if the coordinates are weakly coupled
\end{itemize} \\ \hline
Greedy & $i_k$ is chosen by a greedy rule, such as greatest descent: $$i_k = \argmax_{1 \leq i \leq s} \| [\grad f(\vx^{k - 1})]_{i} \|$$ &
\begin{itemize}[nosep, leftmargin =*]
\item Convergence take the least number of iterations, both theoretically and empirically
\item Well suited for problems with sparse solutions; can be parallelized
\end{itemize} &
\begin{itemize}[nosep, leftmargin =*]
\item Need computation of greedy scores or rankings, such as gradients or difference vectors
\item Highest per-iteration complexity except when greedy scores can be updated at a low cost
\end{itemize} \\ \hline
\end{tabular}
\caption{Summary and comparison of different index choosing schemes for CD.} \label{index rules}
\end{table}
\end{landscape}

\section{Coordinate Friendly Structures} \label{cf structures}

As discussed earlier in the monograph, the strong performance and parallelizability of CD rely on solving subproblems that consist of fewer variables and have low complexities and low memory requirements. Therefore, not all problems are amenable for CD, particularly if little computation is saved from using CD relative to the full update for all coordinates. Intuitively, given $s$ blocks, a block coordinate update should cost about $1/s$ of the computation for a full update. Thus, identifying \textit{coordinate friendly updates} for a given problem is crucial to implementing CD effectively.

In this section, we elaborate on different types of structures in optimization problems that make CD computationally worthy. We define the notion of a \textit{coordinate friendly (CF) update} and \textit{coordinate friendly (CF) structure} and introduce heuristics to exploit problems with CF structures. This basic theory may help practitioners identify when CD is applicable and computationally worthwhile for their problem, as well as determine which quantities to cache and maintain to improve the performance of CD.

The CF structure was originally presented in \cite{peng2016coordinate} for monotone set-valued operators, which apply in more general settings. We refrain from discussing this here and instead replace the notion of an operator with an \textit{update mapping}.

\subsection{Coordinate Friendly Update Mappings} \label{cf mappings}

Before we define coordinate friendly update and structure, we introduce terminology and notation on updates. Suppose we are working with $s$ equally-sized blocks in $\vx\in\RR^n$, i.e. $\vx = (\vx_1, ..., \vx_s)$. Let $T: \RR^{n} \rightarrow \RR^{n}$ represent an \textit{update mapping} or simply \textit{update}. Applying it to $\vx^{k - 1}$, we obtain the next iterate $\vx^k$, i.e.,
$$\vx^k  = T (\vx^{k - 1}).$$
In addition, let $T_{i}$ denote the \textit{coordinate update mapping} of $T$ for block $\vx_{i}$, i.e.,
$$T_i(\vx)=(T(\vx))_i,\,i=1,\ldots,s.$$

As an example, consider the least squares problem:
\begin{equation}\label{ls}
\minimize_{\vx \in \mathbb{R}^{n}} \frac{1}{2}\| \vA \vx - \vb\|_2^2
\end{equation}
where $\vA \in \mathbb{R}^{p \times n}$ and $\vb \in \mathbb{R}^{p}$. The update mapping for gradient descent on (\ref{ls}) is given as:
$$\vx^k=T_{\text{GD}}(\vx^{k-1}) = \vx^{k-1} - \alpha (\vA^T \vA \vx^{k-1} - \vA^T \vb)$$
where $\alpha > 0$ is the step size, 
and the coordinate update mapping is then given as
$$T_{\text{GD},i}(\vx^{k-1})=\vx_i^{k-1}-\alpha (\vA^T \vA \vx^{k-1} - \vA^T \vb)_{i},\,i=1,\ldots,s.$$

We let $\mathcal{N}[a \mapsto b]$ denote the number of basic operations necessary to compute quantity $b$ from the input $a$. Then the update mapping $T$ is called \textit{coordinate friendly (CF)} if
\begin{equation}
\mathcal{N}\big[\vx \mapsto T_{i}(\vx)\big] = O\left(\frac{1}{s} \mathcal{N}[\vx \mapsto T (\vx)]\right),\,\forall i.
\end{equation}
In other words, the number of basic operations necessary to compute the coordinate update is about $1/s$ times of the number of basic operations to compute the full update.

Returning to our least squares example, $T_{\text{GD}}$ is CF since the coordinate gradient update can be computed by
\begin{equation}\label{eq:cgd}T_{\text{GD}, i}(\vx^{k-1})= \vx_{i}^{k - 1} - \alpha \big[(\vA^T \vA)_{i, :} \vx^{k - 1} - (\vA^T \vb)_{i}\big]
\end{equation}
which takes $O(\frac{n^2}{s})$ operations after precomputing $\vA^T \vA$ and $\vA^T \vb$. In contrast, the full gradient update $T_{\text{GD}}$ takes $O(n^2)$ operations with precomputed $\vA^T \vA$ and $\vA^T \vb$.

Intuitively, the definition of coordinate friendly update mappings encapsulates the notion of gaining $s$ times speed-up per coordinate update $T_{i}$ relative to the full update $T$. Equivalently, if we break up our variables into $s$ equally-sized blocks, updating each block should require about $1/s$ times operations necessary to compute the full update.

This formulation also matches our intuition for extreme cases; in the coordinate case of updating individual components $x_{i}$, this will give an $n$ times speed up for computing each coordinate update relative to computing the full update. In the case of only one block containing the entire vector $\vx$, we gain no speed-up.

We say that the problem is \textit{coordinate friendly} or has \textit{coordinate friendly (CF) structure} if there exists a coordinate friendly update for the problem. To readily apply CD, we must recognize and exploit CF structures in optimization problems.

Identifying CF structure in optimization problems is not always trivial. To help with the analysis of optimization problems to find CF structure, we will give some useful heuristics applied to CD implementations in Section \ref{CF heuristics}.

The definitions for CF update and structure may be further generalized for blocks of arbitrary lengths, but we refrain from doing so to maintain simplicity.

\subsection{Common Heuristics for Exploiting CF Structures}\label{CF heuristics}

In this section, we introduce some common heuristics that exploit structure in a problem to improve the performance of CD and gain coordinate friendliness.

\subsubsection{Precomputation of Non-Variable Quantities}

As noted already in the least squares example, one common approach for increasing the computational worthiness of CD is to precompute certain quantities in the update. If certain quantities that do not consist of any of the variables appear in the update, we can precompute them before applying CD to save computation.

For example, consider again the least squares problem given above. 
Since the full gradient is given by $\vA^T \vA \vx - \vA^T \vb$, we can precompute $\vA^T \vA$ and $\vA^T \vb$ to avoid recomputing $\vA^T \vA$ and $\vA^T \vb$ at each iteration.

Note, however, that this is only efficient when $n \sim p$ or $p \gg n$ (recall that $\vA$ has $p$ rows and $n$ columns). If $p \ll n$, then multiplying by $\vA$ then $\vA^T$ is computationally cheaper, so we avoid precomputation in this case.

\subsubsection{Caching and Maintaining Variable-Dependent Quantities}

Another approach to save computation is to cache and maintain variable-dependent quantities. In the CF notation, this would refer to storing some quantity $\mathcal{M}(\vx^k)$ in the memory and updating it at each iteration.

For example, for the least squares problem, instead of performing the coordinate update as in \eqref{eq:cgd}, we can save the quantity
$$\mathcal{M}(\vx^{k - 1}) = \vA^T \vA \vx^{k - 1}.$$
Then for any $i$, $T_{\text{GD},i}$ can be evaluated by
$$T_{\text{GD},i}(\vx^{k-1})=\vx^{k-1}_i-\alpha\left((\mathcal{M}(\vx^{k-1}))_i-(\vA^\top \vb)_i\right),$$
which takes $O(\frac{n}{s})$ operations. Let $\vx^k$ be the vector obtained by block coordinate descent update from $\vx^{k-1}$, i.e.,
\begin{equation}\label{eq:xk-cgd}
\vx^k_i=\left\{
\begin{array}{ll}
T_{\text{GD},i}(\vx^{k-1}), & \text{ if }i=i_k,\\[0.1cm]
\vx_i^{k-1}, & \text{ if }i\neq i_k.
\end{array}\right.
\end{equation}
Then $\mathcal{M}(\vx^k)$ can be obtained by
$$\mathcal{M}(\vx^k)=\mathcal{M}(\vx^{k-1})+(\vA^T \vA)_{:, i_k} (\vx_{i_k}^k - \vx_{i_k}^{k - 1}),$$
which takes $O(\frac{n}{s}) + O(\frac{n^2}{s}) + O(\frac{n}{s}) = O(\frac{n^2}{s})$ operations, since block coordinate addition/subtraction takes $O(\frac{n}{s})$ operations and block matrix-vector multiplication takes $O(\frac{n^2}{s})$ operations. Hence, if $\vA^T \vA$ is cached or $n = O(p)$, we have that
\begin{equation}\label{eq:cf-examp1}\mathcal{N}\big[\{\vx^{k-1}, \mathcal{M}(\vx^{k-1})\}\to \{\vx^k, \mathcal{M}(\vx^k)\}\big]=O\left(\frac{1}{s}\mathcal{N}\big[\vx^{k-1}\to T_{\text{GD}}(\vx^{k-1})\big]\right).
\end{equation}

%

Consider the logistic regression as another example:
$$\minimize_{\vw} F(\vw) = \sum_{j = 1}^m \log \left(1 + \exp[-y_j \vw^T \vx_j]\right), $$
where $\{\vx_j \in \mathbb{R}^n: j=1,\ldots,m\}$ are training data points, and $y_j \in \{-1, 1\}$ is the label of $\vx_j$ for $j = 1, 2, \ldots, m$. The variable $\vw \in \RR^n$ and its block components $\vw_i \in \RR^{n/s}$. Its gradient is
$$\grad F(\vw) = \sum_{j = 1}^m \frac{-y_j \exp[-y_j \vw^T \vx_j]}{1 + \exp[-y_j \vw^T \vx_j]} \vx_j,$$
and the gradient descent update with step size $\alpha$ is given by
\begin{equation}\label{op-lr}
T_{\text{GD}}(\vw)=\vw-\alpha \grad F(\vw).
\end{equation}

To achieve an efficient coordinate update from $\vx^{k-1}$, we maintain the quantity $$\mathcal{M}(\vw^{k - 1}) = \left\{\exp\big[-y_j (\vw^{k-1})^T \vx_j\big],\,j=1,\ldots,m\right\}.$$ Then for any $i$, 
the block gradient descent update from $\vw^{k-1}$ can be computed by
$$T_{\text{GD},i}(\vw^{k-1})=\vw_i^{k-1}-\alpha\sum_{j = 1}^m \frac{-y_j \mathcal{M}(\vw^{k - 1})}{1 + \mathcal{M}(\vw^{k - 1})} (\vx_j)_i,$$
which only takes $O(\frac{mn}{s})$ because computing $\exp[-y_j (\vw^{k - 1})^T \vx_j]$ for $j = 1, ..., m$ has been avoided. Let $\vw^k$ be obtained by applying the block gradient descent update from $\vw^{k-1}$, as given in \eqref{eq:xk-cgd}. Note that
$$\exp\big[-y_j(\vw^k)^\top\vx_j\big]=\exp\big[-y_j(\vw^{k-1})^\top\vx_j\big]\cdot\exp\big[-y_j(\vw^k_{i_k}-\vw^{k-1}_{i_k})^\top(\vx_j)_{i_k}\big],\,\forall j.$$
Thus obtaining $\mathcal{M}(\vw^k)$ from $\mathcal{M}(\vw^{k-1})$ takes $O(\frac{mn}{s})$ operations. Since evaluating $\nabla F(\vw)$ takes $O(mn)$ operations because only the block $\vw_{i_k}$ is needed, we have \eqref{eq:cf-examp1}, and thus $T_{\text{GD}}$ defined in \eqref{op-lr} is CF.

Note that CF structure may also be found in more complicated update mappings, such as the proximal-point update, prox-linear update, and updates derived from primal-dual methods. More details may be found in \cite{peng2016coordinate}.

\section{Applications}\label{applications}

Since we have presented and compared various update schemes and index rules for CD and investigated the coordinate friendliness of problems, we now present canonical examples that apply CD effectively from the literature. These examples demonstrate the efficiency and capability for block CD to handle larger data sets. Each example is initially expressed in its natural form, then re-formulated into a canonical form amenable to block CD methods. A CF analysis is presented for each problem by analyzing the computational cost or number of flops for each coordinate update relative to the full update. Experimental results by the block CD methods are also reported. 

\subsection{LASSO}

In compressed sensing and machine learning, the \textit{Least Absolute Shrinkage and Selection Operator (LASSO)} problem \cite{Tibshirani1996_RegressionShrinkage} seeks to recover a \textit{sparse signal}, or a vector with small support, $\vx^* \in \mathbb{R}^n$ satisfying the underdetermined linear system of equations $\vA \vx = \vb$, where $\vA \in \RR^{m \times n}$ and $\vb \in \RR^m$. It is closely related to the basis pursuit problem \cite{ShaobingChenDonoho1994_BasisPursuit} in signal processing. Under certain conditions, $\ell_1$ minimization returns a sparse solution. This gives the following minimization problem:
\begin{equation*}
\underset{\vx}{\text{minimize}} ~ \|\vx\|_1 + \frac{\lambda}{2}\| \vA \vx - \vb\|_2^2
\end{equation*}
where $\vb \in \mathbb{R}^m$ is our compressed signal, $\lambda \in \mathbb{R}^+$ is a parameter that controls the tradeoff between sparsity and reconstruction, and $\vA \in \mathbb{R}^{m \times n}$ is our measurement matrix with $m \ll n$. Our approach is influenced by the work of \cite{hale2008fixed,li2009coordinate,TsengYun2009_CoordinateGradient}.

\subsubsection{Update Derivation}

Since the problem has a mixed differentiable-nondifferentiable objective, we apply the prox-linear update. Let
\begin{equation*}
F(\vx) = \|\vx\|_1 + \frac{\lambda}{2}\| \vA \vx - \vb\|_2^2
\end{equation*}
Recall that the prox-linear update minimizes the component surrogate function
\begin{equation*}
\vx_{i_k}^k = \argmin_{\vx_{i_k}} f(\vx^{k - 1}) + \langle \grad_{i_k} f(\vx^k), \vx_{i_k} - \vx_{i_k}^{k - 1}\rangle  + \frac{1}{2\alpha_{i_k}^{k - 1}}\|\vx_{i_k} - \vx_{i_k}^{k - 1}\|_2^2 + r_{i_k}(\vx_{i_k}).
\end{equation*}

Note that for the LASSO problem, we let $f(\vx) = \frac{\lambda}{2} \| \vA \vx - \vb\|_2^2$, and thus 
$$\grad f(\vx) = \lambda \vA^T (\vA \vx - \vb).$$
It is also natural to choose our step size as $\alpha_{i_k}^{k - 1} = 1/L_{i_k}$ where $L_{i_k} = \lambda (\vA^T \vA)_{i_k, i_k} = \lambda \|\vA_{:, i_k}\|^2$ is the Lipschitz constant of the gradient of the $i_k$th component of $f(\vx)$. Combined with $r_i(\vx_i) = |\vx_i|$ in the LASSO problem, the prox-linear update becomes
\begin{equation*}
\vx_{i_k}^k = \argmin_{\vx_{i_k}} \|\vA \vx^{k - 1} - \vb\|_2^2 + \langle \lambda \vA^T_{:,i_k} (\vA \vx^{k - 1} - \vb), \vx_{i_k} - \vx_{i_k}^{k - 1}\rangle  + \frac{L_{i_k}}{2}\|\vx_{i_k} - \vx_{i_k}^{k - 1}\|_2^2 + |\vx_{i_k}|.
\end{equation*}

Using the first-order optimality conditions of this minimization, we get that
\begin{align*}
& 0 \in \lambda \vA^T_{:,i_k} (\vA \vx^{k - 1} - \vb) + L_{i_k} (\vx_{i_k} - \vx_{i_k}^{k - 1}) + \partial |\vx_{i_k}|,\\
\mbox{or equivalently,}\quad & \vx_{i_k}^{k - 1} - \frac{\lambda}{L_{i_k}} \vA^T_{:,i_k} (\vA \vx^{k - 1} - \vb) \in (I + \frac{1}{L_{i_k}}\partial | \cdot |)(\vx_{i_k}),
\end{align*}
and thus
\begin{align*}
& \vx_{i_k}^k = \prox_{\frac{1}{L_{i_k}} |\cdot|}(\vx_{i_k}^{k - 1} - \frac{\lambda}{L_{i_k}} \vA^T_{:,i_k} (\vA \vx^{k - 1} - \vb))
\end{align*}
as desired.

Note that $\prox_{\mu |\cdot|}$, which is the proximal operator of the scaled absolute value function $\mu | \cdot |$, is the well-known shrink operator defined by
$$\shrink(x, \mu) = \begin{cases}
x - \mu, & \mbox{ if } x > \mu\\
0, & \mbox{ if } -\mu \leq x \leq \mu\\
x + \mu, & \mbox{ if } x < -\mu
\end{cases}$$
which gives the coordinate update
$$\vx_{i_k}^k = \shrink\big(\vx_{i_k}^{k - 1} - \frac{\lambda}{L_{i_k}} \vA^T_{:,i_k} (\vA \vx^{k - 1} - \vb), \frac{1}{L_{i_k}}\big).$$
Combining this with $L_{i_k} = \lambda \| \vA_{:, i_k}\|^2$ gives the final coordinate update
$$\vx_{i_k}^k = \shrink\big(\vx_{i_k}^{k - 1} - \frac{1}{\|\vA_{:, i_k}\|^2} \vA^T_{:,i_k} (\vA \vx^{k - 1} - \vb), \frac{1}{\lambda \|\vA_{:, i_k}\|^2}\big).$$

\subsubsection{Continuation}

Since the step size $1/L$ may be small, we may improve the speed of convergence of the algorithm by introducing \textit{continuation}. Note that smaller values of $\lambda$ dictate sparser solutions while larger values of $\lambda$ admit less sparse solutions. This remark motivates defining a series of problems with increasing $\lambda_k$ that reaches the final $\lambda$. Each problem is solved consecutively, with the solution to the current problem acting as the starting point for the next problem. A simple algorithm for $\lambda_{k + 1}$ is given by
$$\lambda_{k + 1} = \eta \lambda_k$$
where $\eta > 1$. For more details, please refer to \cite{hale2008fixed}.

\subsubsection{Derivations for Gauss-Southwell Rules}

We also show the derivation for the Gauss-Southwell rules:

\begin{enumerate}
\item \textit{GS-s index rule}: Recall that the GS-s rule chooses an index $i_k$ at each iteration $k$ by (\ref{gs-s}). Note that for the LASSO problem,
$$\grad_j f(\vx^{k - 1}) = \lambda \vA^T_{:,j} (\vA \vx - \vb)$$
$$\partial r_j(\vx^{k - 1}_{j}) = \begin{cases}
1 & \mbox{if } \vx_j^{k - 1} > 0\\
[-1, 1] & \mbox{if } \vx_j^{k - 1} = 0\\
-1 & \mbox{if } \vx_j^{k - 1} < 0.
\end{cases}$$
Therefore, following $g_j(\vx^{k - 1}) = \min_{\tilde{\grad} r_j \in \partial r_j} \| \grad_j f(\vx^{k - 1}) + \tilde{\grad} r_j (\vx_j^{k - 1}) \|$, we get
$$g_j(\vx^{k - 1}) = \begin{cases}
\| \lambda \vA^T_{:,j} (\vA \vx^{k - 1} - \vb) + \sign(\vx_j^{k - 1})\| & \mbox{if } \vx_j^{k - 1} \neq 0\\
\| \shrink(\lambda \vA^T_{:,j} (\vA \vx^{k - 1} - \vb), 1) \| & \mbox{otherwise.}
\end{cases}$$
Choosing the largest score $g_j (\vx^{k - 1})$ gives the index $i_k$.

\item \textit{GS-r index rule}: Recall that the GS-r rule chooses an index $i_k$ at each iteration $k$ by (\ref{gs-r}). Therefore,
$$i_k = \argmax_j \big\|\vx_j^{k - 1} - \prox_{\frac{1}{L}|\cdot|} (\vx_j^{k - 1} - \frac{\lambda}{L} \vA^T_{:,j} (\vA \vx^{k - 1} - \vb))\big\|$$
where $L = \lambda \|\vA^T \vA\|$ is the Lipschitz constant of the smooth quadratic function.

\item \textit{GS-q index rule}: Recall that the GS-q rule chooses an index $i_k$ at each iteration $k$ by (\ref{gs-q}). Since
$$d_j = \prox_{\frac{1}{L}|\cdot|} (\vx_j^{k - 1} - \frac{\lambda}{L} \vA^T_{:,j} (\vA \vx^{k - 1} - \vb)) - \vx_j^{k - 1}$$
for the LASSO problem, we can plug in $\vd$ into our equation
$$\frac{\lambda}{2} \| \vA \vx^{k - 1} - \vb\|_2^2 + \lambda d_j (\vA^T_{:,j} (\vA \vx^{k - 1} - \vb))  + \frac{L}{2} |d_j|^2 + |\vx_j^{k - 1} + d_j| - |\vx_j^{k - 1}|.$$
Finding the index for the smallest score gives the index $i_k$. Note that the first term is constant in $j$, so it can be dropped from the score computation.
\end{enumerate}

\subsubsection{CF Analysis}

Following our analysis of least squares in Section \ref{cf mappings}, we give a CF analysis of the LASSO problem. Note that though the LASSO and least squares problems are similar, $\vA \in \RR^{m \times n}$ is a short-and-fat matrix, i.e. $m \ll n$. Therefore, the precomputation of $\vA^T \vA$ and $\vA^T \vb$ is not incentivized.

Instead, we cache and maintain the quantity $\vA \vx^k$. Recall that for the full update, we can simply multiply first by $\vA$ then $\vA^T$ to take advantage of our short-and-fat matrix $\vA$:
$$\vx^k = \shrink(\vx^{k - 1} - \frac{1}{L} \vA^T (\vA \vx^{k - 1} - \vb), \frac{1}{\lambda L}).$$
This gives $\mathcal{N}[\vx^k \mapsto T(\vx^k) ] = O(mn)$ since matrix-vector multiplication and vector addition takes $O(mn)$ operations and $O(m)$ or $O(n)$ operations, respectively.

For the coordinate update, caching  $\vA \vx^{k-1}$ gives the cheap coordinate update
$$\vx_{i_k}^k = \shrink(\vx_{i_k}^{k - 1} - \frac{1}{L_{i_k}} (\vA_{:,i_k}^T (\vA \vx^{k - 1} - \vb)), \frac{1}{\lambda L_{i_k}}),$$
where computing $\vA_{:,i_k}^T (\vA \vx^{k - 1} - \vb)$ involves only two $O(m)$ vector-vector operations.
We maintain $\vA \vx^k$ by adding $(\vx_{i_k}^k - \vx_{i_k}^{k - 1}) \vA_{:, i_k}$ to $\vA \vx^{k - 1}$, i.e.
$$\vA \vx^k = \vA \vx^{k - 1} + (\vx_{i_k}^k - \vx_{i_k}^{k - 1}) \vA_{:, i_k}.$$
Since both the coordinate update and maintenance of $\vA \vx^k$ take $O(m)$ operations,
$$\mathcal{N}[\{\vx^k, \vA \vx^k\} \mapsto \{T_{i_k}(\vx^k), \vA T_{i_k}(\vx^k) \} ] = O(m) = O(\frac{1}{n} \mathcal{N}[\vx^k \mapsto T(\vx^k) ])$$
which shows that this approach is CF.

\subsubsection{Numerical Examples}

For the following numerical examples, we apply our approach to a compressed, randomly generated, sparse signal, and compare the performance of the algorithm for multiple index rules. We first specify the size and type of matrix $\vA \in \RR^{m \times n}$, $k = |\supp (\vx_s)|$, and standard deviation of noise $\sigma$. The entries of the matrix $\vA$ are taken from a standard normal distribution. We choose our support $\supp (\vx_s)$ by using a random permutation, and we set non-zero entries following a normal distribution $\mathcal{N}(0, 2)$.

Next, we introduce Gaussian noise to our compressed signal $\vb = \vA \vx_s + \epsilon$, where $\epsilon \sim \mathcal{N}(0, \sigma)$ is a Gaussian noise vector.

\begin{figure}
\center
\includegraphics[scale = 0.5]{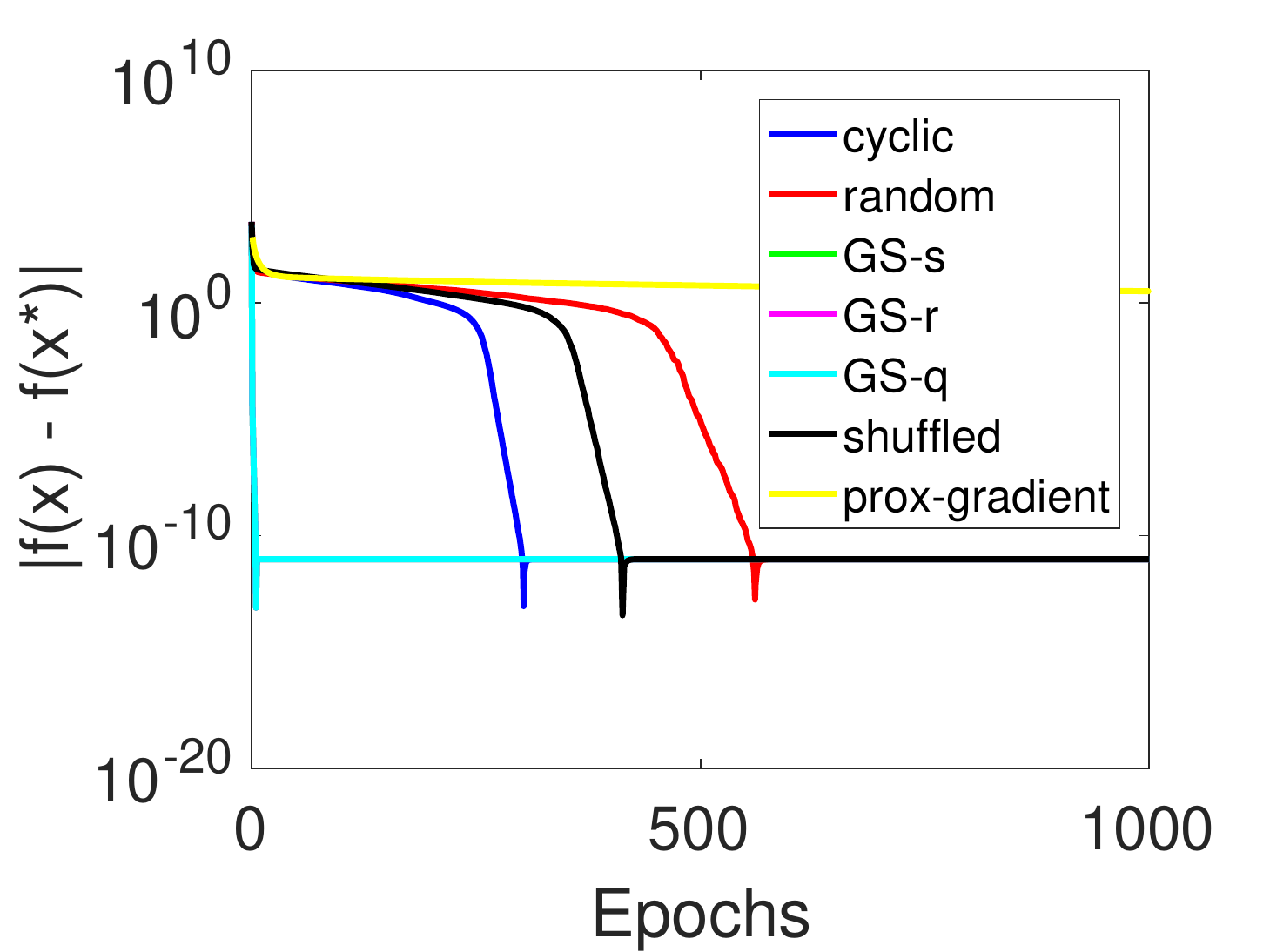}
\includegraphics[scale = 0.5]{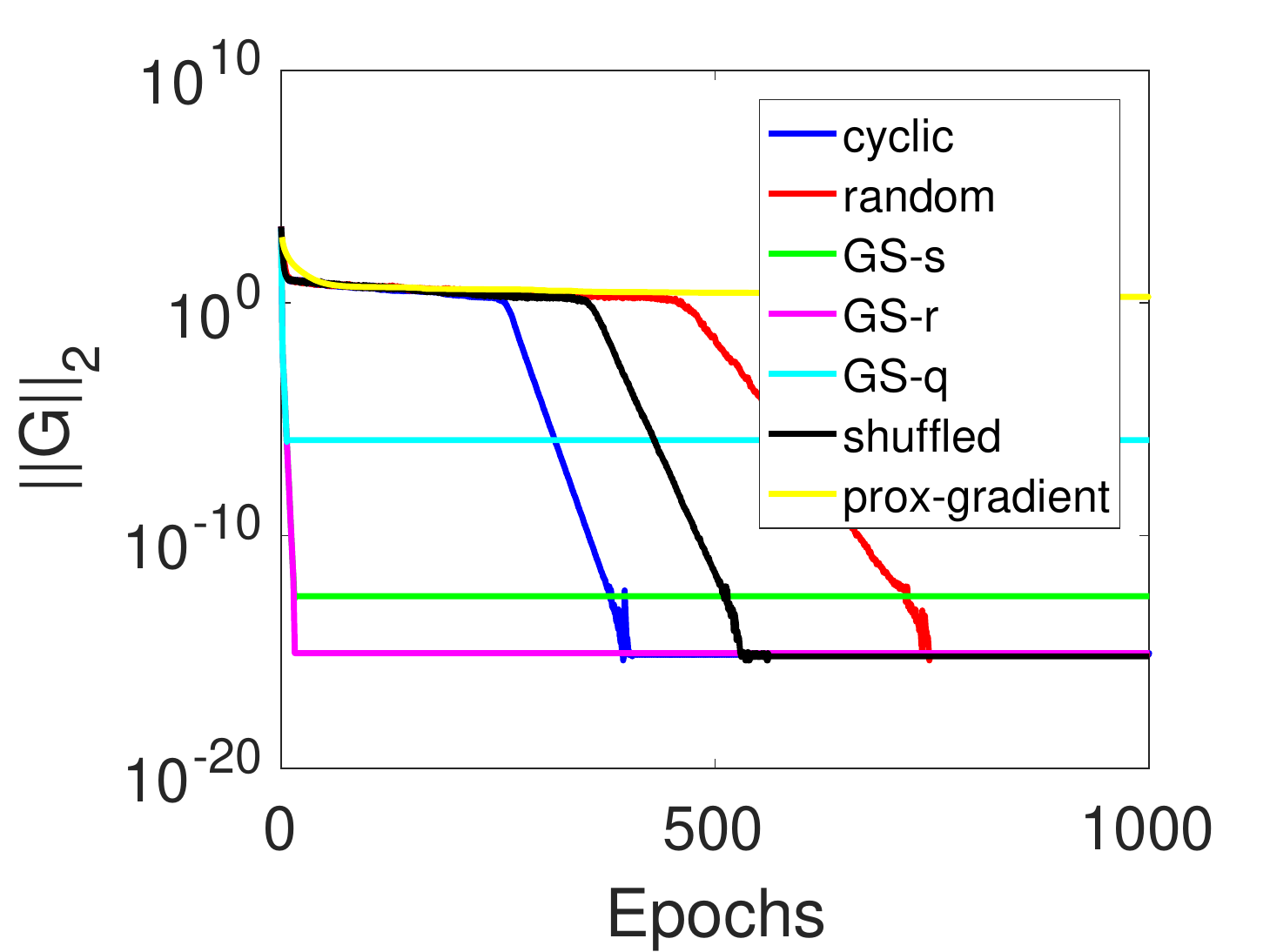}
\includegraphics[scale = 0.5]{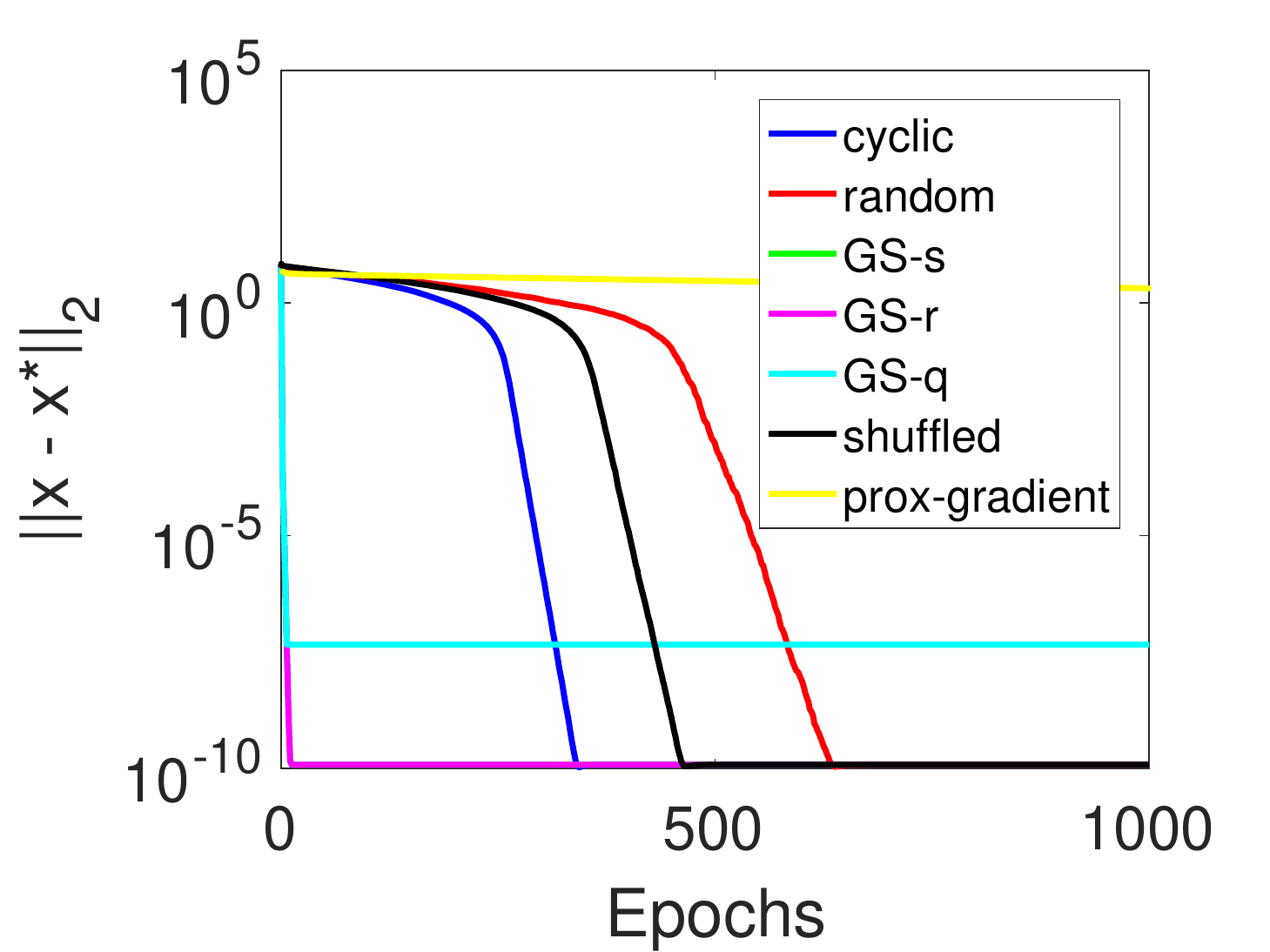}
\caption{Convergence results for LASSO application with $m = 50$, $n = 100$, and $\sigma = 10^{-4}$. Randomized and shuffled cyclic CD results are averaged over 100 trials. The top left, top right, and bottom graphs compare the objective value $\|\vx\|_1 + \frac{\lambda}{2}\|\vA \vx - \vb\|_2^2$, norm of the gradient map, and distance to the minimizer against the number of epochs, or $n$ coordinate updates.}\label{fig: lasso}
\end{figure}

For our experiments, we set $m = 50$ and $n = 100$, and use $\sigma = 10^{-4}$ and $\lambda = 10^3$. We plot both the objective, norm of the gradient map, and distance to the solution in Figure \ref{fig: lasso}. Note that the gradient map here is defined as
$$\vG = \vx - \prox_{\frac{1}{L}\| \cdot \|_1}(\vx - \frac{1}{L} \grad f(\vx)).$$
We average the results of randomized and shuffled cyclic CD over 100 trials. We use CVX with high precision to find $\vx^*$ for comparison \cite{grant2008cvx}.

Due to the sparsity of the solutions, the Gauss-Southwell rules perform significantly better than their randomized, shuffled cyclic, and cyclic counterparts. Cyclic and shuffled cyclic variants also perform better than randomized coordinate descent for this particular application, since randomized coordinate descent does not discover the support of the solution. The drop in objective also suggests that greedy and cyclic CD uses the first number of iterations to discover the support, then solves for the solution over the support, whereas randomized CD fails to exploit the sparsity of the solution.

\subsection{Non-Negative Matrix Factorization for Data Mining and Dimensionality Reduction}

Non-negative Matrix Factorization (NMF) \cite{paatero1994positive, cichocki2009nonnegative} aims at factorizing a non-negative matrix $\vM \in \RR_+^{m \times n}$ into two non-negative matrices of column-size $r$, i.e. finding non-negative $\vX \in \RR_+^{m \times r}$ and $\vY \in \RR_+^{n \times r}$ such that
\begin{equation}\label{eq:nmf-goal}
\vM = \vX\vY^T, \text{ or }\vM\approx \vX\vY^T.
\end{equation}
When $r$ is significantly smaller than $\min\{m,n\}$, the matrix $\vM$  is (approximately) low rank. Since the resulting matrices $\vX$ and $\vY$ are both non-negative, they interpret meaningfully and are easy to inspect. This particular advantage motivates the use of NMF in a variety of fields, including text mining, computer vision, recommender systems, and signal processing.

The goal in \eqref{eq:nmf-goal} can be achieved by solving the following minimization problem:
\begin{equation}\label{eq:nmf-ls}
\begin{aligned}
& \underset{\vX, \vY}{\text{minimize}}
& & \frac{1}{2} \|\vM - \vX \vY^T \|_F^2 \\
& \text{subject to}
& & \vX, \vY \geq 0
\end{aligned}
\end{equation}
Note that \eqref{eq:nmf-ls} is a non-convex problem due to the bilinear term $\vX\vY^T$. Toward finding a solution to \eqref{eq:nmf-ls}, we follow  \cite{xu2013block, peng2016coordinate} and apply block coordinate descent with prox-linear update. By incorporating the constraints into the objective function as indicator variables, we obtain the problem
\begin{equation}\label{eq:nmf-ls-ind}
\minimize_{\vX, \vY} \frac{1}{2} \| \vM - \vX \vY^T \|_F^2 + \iota_{\geq 0}(\vX) + \iota_{\geq 0}(\vY),
\end{equation}
where
$$\iota_{\geq 0}(\vX) = \begin{cases} 0 & \mbox{ if } \vX_{ij} \geq 0 \mbox{ for all } i,j \\ \infty & \mbox{otherwise.} \end{cases}$$

\subsubsection{Update Derivation}


Let
$$F(\vX, \vY) = \frac{1}{2}\| \vM - \vX \vY^T \|_F^2.$$
The gradient can be easily obtained as
\begin{equation}\label{eq:nmf-full-grad}
\nabla F(\vX,\vY)=\left[\nabla_\vX F(\vX,\vY),\nabla_\vY F(\vX,\vY)\right]=\left[(\vX \vY^T - \vM)\vY,\,(\vX \vY^T - \vM)^T \vX\right].
\end{equation}
To enforce nonnegativity, the projected gradient method goes a step along the gradient descent direction and then projects the iterate to the nonnegative orthant, i.e.,
\begin{subequations}\label{eq:nmf-pj-grad}
\begin{align}
&\vX^k = \max\left(0, \vX^{k - 1} - \eta_{k - 1} \grad_{\vX} F(\vX^{k - 1}, \vY^{k - 1})\right),\\
&\vY^k = \max\left(0, \vY^{k - 1} - \eta_{k - 1} \grad_{\vY} F(\vX^{k - 1}, \vY^{k - 1})\right),
\end{align}
\end{subequations}
where $\eta_{k - 1}$ is the step size that can be determined by line search. 

To derive the projected block coordinate gradient update, we partition the variables into disjoint blocks by columns, i.e. $\vX = (\vX_{:, 1}, \ldots, \vX_{:, r})$ and $\vY = (\vY_{:, 1}, \ldots, \vY_{:, r})$. (Partition by rows or by block is also possible.) At each iteration $k$, we select one block variable and perform one of the following updates: 
\begin{subequations}\label{eq:nmf-pj-pgrad}
\begin{align}
&\vX_{:, i_k}^k = \max\left(0, \vX_{:, i_k}^{k - 1} - \eta_{k - 1} \grad_{\vX_{:, i_k}} F(\vX^{k - 1}, \vY^{k - 1})\right),\text{ if }i_k\text{th column of }\vX\text{ selected}\\
&\vY_{:, i_k}^k = \max\left(0, \vY_{:, i_k}^{k - 1} - \eta_{k - 1} \grad_{\vY_{:, i_k}} F(\vX^{k - 1}, \vY^{k - 1})\right), \text{ if }i_k\text{th column of }\vY\text{ selected},
\end{align}
\end{subequations}
where the partial gradients can be explicitly evaluated via
\begin{subequations}\label{eq:nmf-partial-grad}
\begin{align}
&\grad_{\vX_{:, i_k}} F(\vX, \vY) = \vX (\vY^T\vY_{:, i_k}) - \vM\vY_{:, i_k},\\
&\grad_{\vY_{:, i_k}} F(\vX, \vY) = \vY (\vX^T\vX_{:, i_k}) - \vM^T \vX_{:, i_k},
\end{align}
\end{subequations}
and the step size $\eta_{k-1}$ can be set to the reciprocal of the Lipschitz constant of the corresponding partial gradient, i.e.,
\begin{equation}\label{eq:step-size-eta}\eta_{k-1}=\left\{
\begin{array}{ll}
\frac{1}{\|\vY_{:, i_k}^{k - 1}\|_2^2},&\text{ for }\vX\text{-update},\\[0.4cm]
\frac{1}{\|\vX_{:, i_k}^{k - 1}\|_2^2},&\text{ for }\vY\text{-update}.
\end{array}\right.
\end{equation}
To avoid \emph{zero}-columns during each iteration that would give an infinite step sizes in \eqref{eq:step-size-eta}, we can 
simply restrict $\vX$ to have unit-norm columns, since $\vX \vY^T = (\vX \vD)(\vY \vD^{-1})^T$ for any invertible diagonal matrix $\vD$; see \cite{xu2014globally} for more details.

Note that 
 if the subproblems in \eqref{eq:nmf-pj-pgrad} are still expensive for extremely large problems, CD may further break $\vX$ and $\vY$ into more blocks that can be updated in a sequential fashion. One can also apply direct minimization of the objective rather than the prox-linear approach described above. That causes each subproblem to be a non-negative least squares problem, which has many off-the-shelf solvers.

\subsubsection{Derivations for Gauss-Southwell Rules}

We give the derivations for the Gauss-Southwell rules:

\begin{enumerate}
\item \textit{GS-s index rule}: Recall that the GS-s rule chooses an index $i_k$ at each iteration $k$ by (\ref{gs-s}). Since the subdifferential of the indicator function of a closed convex set $C$ is the normal cone, i.e.
$$\partial \iota_C (\vx) = N_C (\vx) = \{\vg \in \mathbb{R}^n : \vg^T \vx \geq \vg^T \vy \mbox{ for all } \vy \in C\},$$
the subdifferential of $\iota_{\geq 0}$ at any $\vx\ge0$ is given as
$$\partial \iota_{\geq 0} (\vx)=\{\vg\le 0: g_ix_i=0,\forall i\}.$$
Since $\vX^k,\vY^k\ge0,\forall k$, the indices by GS-s rule for $\vX$ and $\vY$ are given respectively by
$$i_k = \argmax_j \| (\vG_{\vX})_{:,j}\| $$
$$i_k = \argmax_j \| (\vG_{\vY})_{:,j}\|$$
where
$$(\vG_{\vX})_{i, j} = \begin{cases}
\min (\grad_{\vX_{i,j}}F (\vX^{k - 1}, \vY^{k - 1}), 0) & \mbox{ if } \vX_{i, j}^{k-1} = 0\\
\grad_{\vX_{i,j}} F(\vX^{k - 1}, \vY^{k - 1}) & \mbox{ otherwise}
\end{cases}$$
$$(\vG_{\vY})_{i, j} = \begin{cases}
\min (\grad_{\vY_{i,j}}F (\vX^{k - 1}, \vY^{k - 1}), 0) & \mbox{ if } \vY_{i, j}^{k-1} = 0\\
\grad_{\vY_{i,j}} F(\vX^{k - 1}, \vY^{k - 1}) & \mbox{ otherwise.}
\end{cases}$$
\item \textit{GS-r index rule}: Recall that the GS-r rule chooses an index $i_k$ at each iteration $k$ by (\ref{gs-r}). Thus,
$$i_k = \argmax_j \| \vX_{:, j}^{k - 1} - \max\left(0, \vX_{:, j}^{k - 1} - \eta_{k - 1} \grad_{\vX_{:, j}} F(\vX^{k - 1}, \vY^{k - 1})\right)\| $$
$$i_k = \argmax_j \| \vY_{:, j}^{k - 1} - \max\left(0, \vY_{:, j}^{k - 1} - \eta_{k - 1} \grad_{\vY_{:, j}} F(\vX^{k - 1}, \vY^{k - 1})\right)\| $$
are the selected indices for $\vX$ and $\vY$, respectively.
\item \textit{GS-q index rule}: Recall that the GS-q rule chooses an index $i_k$ at each iteration $k$ by (\ref{gs-q}). Let
\begin{align*}
&\vd_j^x=\max\left(0, \vX_{:, j}^{k - 1} - \eta_{k - 1} \grad_{\vX_{:, j}} F(\vX^{k - 1}, \vY^{k - 1})\right) - \vX_{:, j}^{k - 1},\,j=1,\ldots,r,\\
&\vd_j^y=\max\left(0,  \vY_{:, j}^{k - 1} - \eta_{k - 1} \grad_{\vY_{:, j}} F(\vX^{k - 1}, \vY^{k - 1})\right) - \vY_{:, j}^{k - 1},\,j=1,\ldots,r.
\end{align*}
Then the indices for $\vX$ and $\vY$ are given respectively by
\begin{align*}
&i_k=\argmin_j \grad_{\vX_{:, j}} F(\vX^{k - 1}, \vY^{k - 1})^T \vd_j^x + \frac{\eta_{k-1}}{2} \| \vd_j^x\|^2,\\
&i_k=\argmin_j \grad_{\vY_{:, j}} F(\vX^{k - 1}, \vY^{k - 1})^T \vd_j^y + \frac{\eta_{k-1}}{2} \| \vd_j^y\|^2,
\end{align*}
where we have used the fact that $\vX_j^{k - 1} + \vd_j^x$ and $\vY_j^{k - 1} + \vd_j^y$ are both componentwise non-negative.
\end{enumerate}

\subsubsection{CF Analysis}

To show the CF property of the projected gradient method, we compare its per-update cost to that of the projected coordinate gradient method.  
Note that the full gradient update in \eqref{eq:nmf-full-grad} 
costs $O(mnr)$ operations, and with the full gradient computed, the nonnegativity projection in \eqref{eq:nmf-pj-grad} takes another $O(mn)$ operations. Hence, one update of the projected gradient costs $O(mnr)$. 

For the projected coordinate gradient update, 
the evaluation of partial gradient in \eqref{eq:nmf-partial-grad} costs $O(mn)$, and with the partial gradient computed, the projection in \eqref{eq:nmf-pj-pgrad} takes another $O(m)$ for the $\vX$-update and $O(n)$ for the $\vY$-update. 
Hence, one update of the projected coordinate gradient costs $O(mn)$, and according to our definition of CF property in section \ref{cf structures}, the projected gradient mapping is CF.


\subsubsection{Numerical Example}

\begin{figure}
\center
\includegraphics[scale=.5]{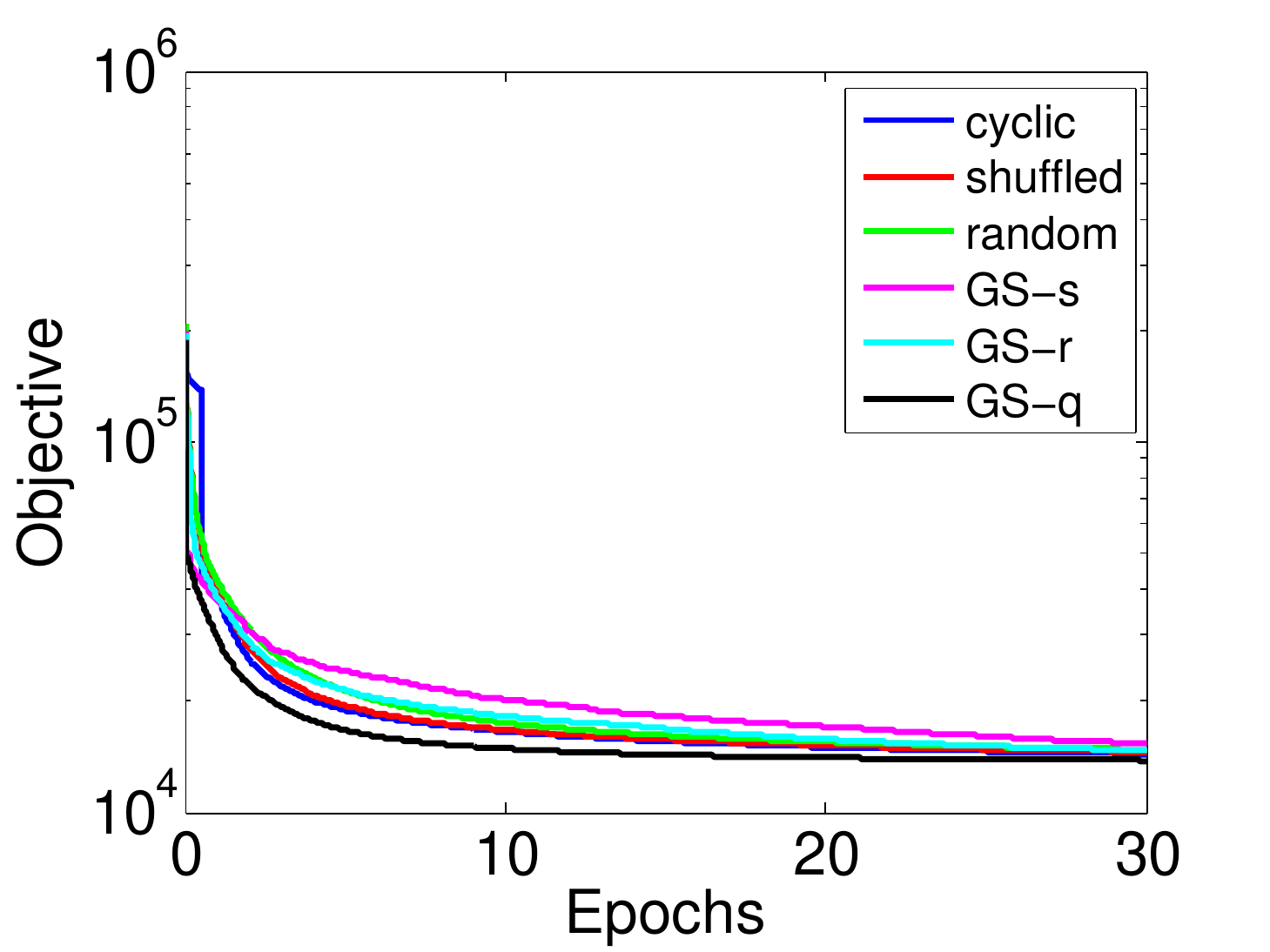}
\includegraphics[scale=.5]{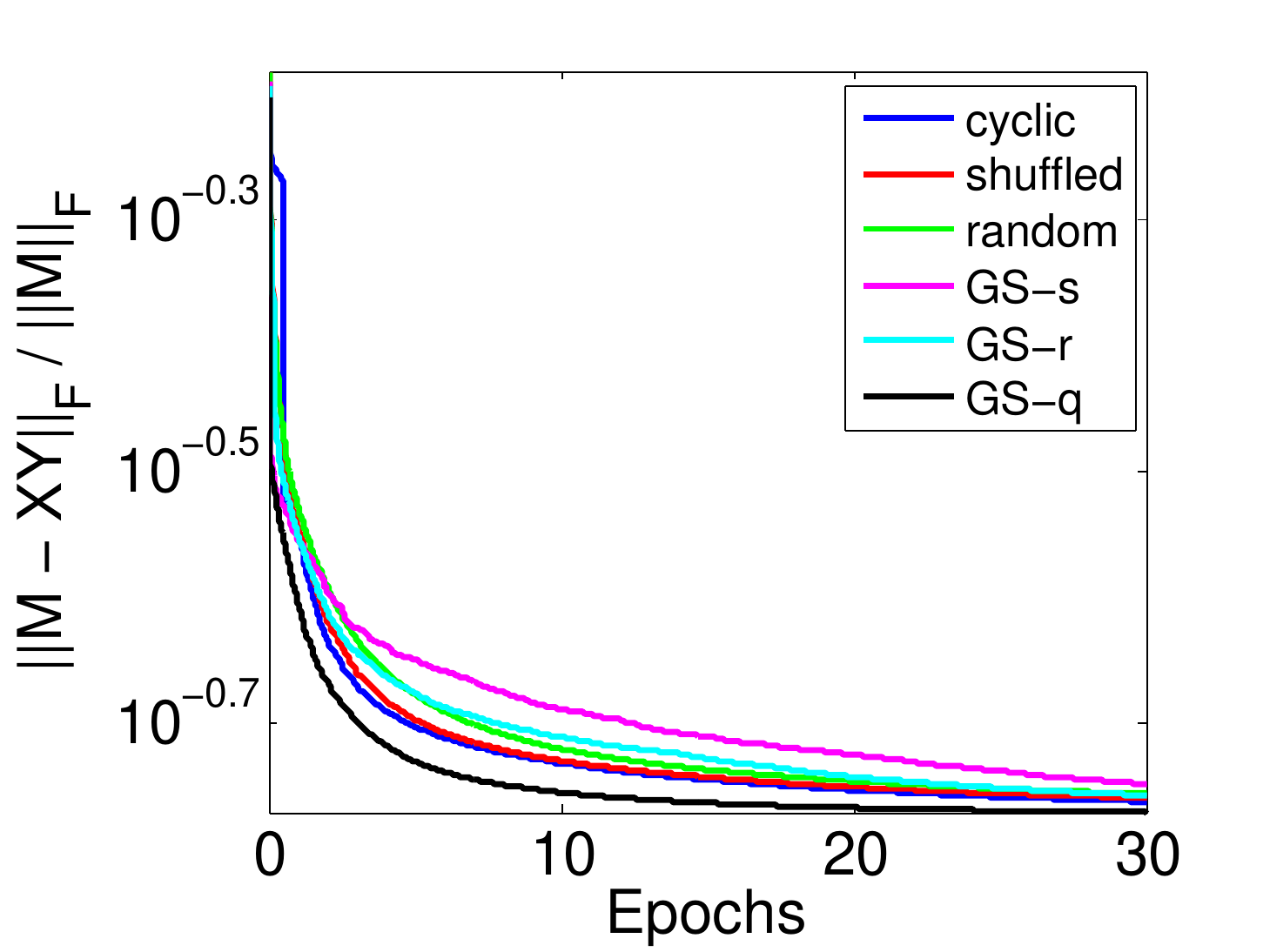}
\caption{Convergence results comparing cyclic, shuffled cyclic, random, and greedy projected coordinate descent for the NMF problem \eqref{eq:nmf-ls} on the ORL data set. The graphs compare the objective value $\frac{1}{2}\| \vM - \vX \vY^T\|_F^2$ and relative error $\|\vM - \vX \vY^T\|_F^2/\|\vM\|_F$ against the number of epochs.}\label{nmf exp}
\end{figure}

We apply the projected coordinate gradient update in \eqref{eq:nmf-pj-pgrad} to the NMF problem \eqref{eq:nmf-ls} on the ORL database from AT\&T Laboratories Cambridge. This data set consists of 40 different faces taken from 10 different directions and with different expressions. We vectorize each image and obtain a matrix $\vM$ of size $10304 \times 400$. We test its performance using cyclic, shuffled cyclic, random, and all greedy index rules. For the performance of randomized projected CD, we average the results from 100 trials. The objective and relative error decrease are plotted in Figure \ref{nmf exp}. 


The GS-q outperforms all other rules while the shuffled cyclic and cyclic rules edge out the other rules. GS-s and GS-r perform surprisingly poor here despite their complexity and are comparable to randomized CD. This observation may suggest that shuffled cyclic or cyclic rules may be better than more expensive Gauss-Southwell rules for NMF.

\subsection{Sparse Logistic Regression for Classification}\label{LR}
In machine learning and statistics, we are often concerned with predicting the categories of new observations. One common approach to do this is to learn a model from a training data set with correctly classified observations, and then use the learned model to predict category membership for new observations.

Logistic regression \cite{HosmerJr.LemeshowSturdivant2013_IntroductionLogistic} is one popularly used model. 
It estimates the probability of a binary response based on given multivariable data points. 
To fit the model, one can solve an optimization problem to search its parameters. We describe our approach below. 

Given a set of training data
\begin{equation} \label{eq:training}
\mathcal{D} = \{(\bx_i,y_i)| ~\bx_i \in \R^n,~y_i \in \{-1,+1\} \},
\end{equation}
consisting of features $\vx_i$ and labels $y_i$, $\ell_1$-regularized logistic regression \cite{Tibshirani1996_RegressionShrinkage} can be formulated as the following optimization problem:
\begin{equation*} \label{eq:LR}
\begin{aligned}
& \minimize_{\vw} &\| \vw \|_1 + C \sum \limits_{i=1}^{l}  \log (1+e^{-y_i \vw^T \bx_i})\\
\end{aligned}
\end{equation*}
where $C>0$ is a parameter that controls the balance between regularization and loss function. Note that the $\ell_1$-regularization term can be replaced by some other regularizer depending on the application. 

\subsubsection{Update Derivation}
We follow \cite{WuLange2008_CoordinateDescent,Logistic} and apply the \textit{one-dimensional Newton direction} coordinate update to solve the sparse logistic regression problem \eqref{eq:LR}.
Let
\begin{equation*}
F(\vw) = \| \vw \|_1 + C \sum \limits_{i=1}^{l}  \log (1+e^{-y_i \vw^T \bx_i})
\end{equation*}
Instead of updating $\vw$, we update one coordinate $\vw_{i_k}$. Therefore, we solve the one-variable subproblem:
\begin{equation*} \label{minLR}
\begin{split}
\minimize_d ~ g_{i_k}(d) & = F(\vw^k+d \be_{i_k}) -F(\vw^k) \\
& = | \vw^k_{i_k}+d|-|\vw^k_{i_k}| + f(\vw^k+d \be_{i_k}) - f(\vw^k)
\end{split}
\end{equation*}
where $$f(\vw)=C \sum \limits_{i=1}^{l}  \log (1+e^{-y_i \vw^T \bx_i}).$$ 

At the $k$th iteration, given index $i_k$, we first find a coordinate descent direction by minimizing the second-order approximation of the smooth function,
\begin{equation} \label{eq:CDN}
\minimize_d | \vw^k_{i_k}+d|-|\vw_{i_k}| + f'_{i_k}(\vw^k) d + \frac{1}{2} f''_{i_k}(\vw^k) d^2
\end{equation}
where
$f_{i_k}'(\vw)$ and $f_{i_k}''(\vw)$ denote the first and second derivative with respect to the $i_k$th component, respectively. It is easy to see that \eqref{eq:CDN} has a closed-form solution by using the shrinkage operator:

\[ d = \left\{
  \begin{array}{l l}
    -\frac{f'_{i_k}(\vw^k)+1}{f''_{i_k}(\vw^k)} & \quad \text{if $f'_{i_k}(\vw^k)+1 \leq f''_{i_k}(\vw^k) \vw_{i_k}^k$}\\
     -\frac{f'_{i_k}(\vw^k)-1}{f''_{i_k}(\vw^k)} & \quad \text{if $f'_{i_k}(\vw^k)-1 \geq f''_{i_k}(\vw^k) \vw_{i_k}^k$}\\
-\vw_{i_k}^k & \quad \text{otherwise.}
  \end{array} \right. \]

To guarantee a decreasing in the function value, armijo rule \cite{tseng2009coordinate} is applied to find a stepsize $\alpha^k \in (0,1)$ such that:
\begin{equation*}
g_{i_k}(\alpha^k d) \leq \sigma \alpha^k (f'_{i_k}(\vw^k)d + |\vw^{k}_{i_k}+d|-|\vw^k_{i_k}|)
\end{equation*}
where $\sigma$ is an arbitrary constant in (0,1).
Thus we get the update
\begin{equation} \label{eq:Logistic_update}
\vw^{k+1}_{i_k} = \vw^k_{i_k}+\alpha^k d.
\end{equation}
Note that $\vw_{i_k}^k$ reaches the optimal value when  $\frac{\partial g_{i_k}}{\partial d}(0) = 0$.\\

\subsubsection{Derivations for Gauss-Southwell Rules}
We derive the update schemes for the Gauss-Southwell rules:
\begin{enumerate}
\item \textit{GS-s index rule}: Recall that the GS-s rule chooses an index $i_k$ at each iteration $k$ by (\ref{gs-s}). Note for the Logistic Regression problem, we have $f(\vw) = C\sum_{i=1}^l \log (1+e^{-y_i \vw^T \bx_i})$ and $r_i(\vw) = |\vw_i|$, and thus
$$ \nabla_j f(\vw^{k-1}) = C\sum_{i=1}^l \frac{-y_i x_{ij} e^{-y_i \vw^T \bx_i}}{1+e^{-y_i \vw^T \bx_i}},  $$
\[ \partial r_j(\vw_j^{k-1}) =
\begin{cases}
1 &\quad \text{ if } \vw_j > 0 \\
[-1,1] &\quad \text{ if } \vw_j = 0 \\
-1 &\quad \text{ if } \vw_j < 0, \\
\end{cases}
\]
where $x_{ij}$ is the $j$th entry of $\bx_i$.
Thus, if $g_j(\vw^{k-1}) = \min_{\tilde{\nabla}{r_j}\in \partial r_j} \| \nabla_j f(\vw^{k-1}) + \tilde{\nabla} r_j (\vw^{k-1}_j)\|$,
\[ g_j(\vw^{k-1}) =
\begin{cases}
\| C\sum_{i=1}^l \frac{-y_i x_{ij}e^{-y_i \vw^T \bx_i}}{1+e^{-y_i \vw^T \bx_i}} + \text{ sign}(\vw_j^{k-1})\| &\quad \text{ if } \vw_j \neq 0 \\
\| \text{ shrink } ( C\sum_{i=1}^l \frac{-y_i x_{ij}e^{-y_i \vw^T \bx_i}}{1+e^{-y_i \vw^T \bx_i}},1)\| &\quad \text{ otherwise.}
\end{cases}
\]
Choosing the largest score $|g_j(\vw^{k-1})|$ gives the index $i_k$.
\item \textit{GS-r index rule}:  Recall that the GS-r rule chooses an index $i_k$ at each iteration $k$ by (\ref{gs-r}). Thus,
$$ i_k = \argmax_j \|\vw_j^{k-1}  - \prox_{\frac{1}{L}|\cdot |}(\vw_j^{k-1}-\frac{C}{L} \sum_{i=1}^l \frac{-y_i x_{ij}e^{-y_i \vw^T \bx_i}}{1+e^{-y_i \vw^T \bx_i}})\|,$$
where $L$ is the Lipschitz constant of $\nabla f$.
Choosing the largest score $|d|$ gives the index $i_k$.
\item \textit{GS-q index rule}: Recall that the GS-q rule chooses an index $i_k$ at each iteration $k$ by (\ref{gs-q}). Since
$$ d= \prox_{\frac{1}{L}|\cdot |}(\vw_j^{k-1}-\frac{C}{L} \sum_{i=1}^l \frac{-y_i x_{ij}e^{-y_i \vw^T \bx_i}}{1+e^{-y_i \vw^T \bx_i}})-\vw_j^{k-1}, $$
we can plug $d$ into equation
$$ f(\vw^{k-1}) +d \nabla_j f(\vw^{k-1})  + \frac{L}{2} |d|^2 + |\vw_j^{k-1}+d| - |\vw_j^{k-1}|$$
Finding the index for the smallest score gives the index $i_k$.
\end{enumerate}

\subsubsection{CF analysis}
Note that
\begin{equation*}
\begin{split}
f'_{j}(\vw^{k-1}) & = C\sum_{i=1}^l y_i  x_{ij}(\frac{1}{1+e^{-y_i\vw^{k-1}\bx_i}}-1), \\
f''_{j}(\vw^{k-1}) & = C\sum_{i=1}^l x_{ij}^2 (\frac{1}{1+e^{-y_i\vw^{k-1}\bx_i}})(1-\frac{1}{1+e^{-y_i\vw^{k-1}\bx_i}}),
\end{split}
\end{equation*}
where we have eliminated $y_i^2\equiv 1$.
We cache $e^{-y_i\vw^{k-1}\bx_i}$ and, following \eqref{eq:Logistic_update}, update $e^{-y_i\vw^{k}\bx_i}$ by
\begin{equation*}
e^{-y_i\vw^{k}\bx_i} = e^{-y_i\vw^{k-1}\bx_i}\cdot e^{-y_i\lambda d x_{ij}}
\end{equation*}
Then the computations take $O(l)$ operations, so
$$\mathcal{N}[\{\vw^k_{i_k} \mapsto \{T_{i_k}(\vw^k_{i_k}) \} ] = O(l) = O(\frac{1}{l} \mathcal{N}[\vw^k \mapsto T(\vw^k) ]),$$
which shows that this problem has CF structure.

\subsubsection{Numerical Example}
\begin{figure}
\center
\includegraphics[scale=0.6]{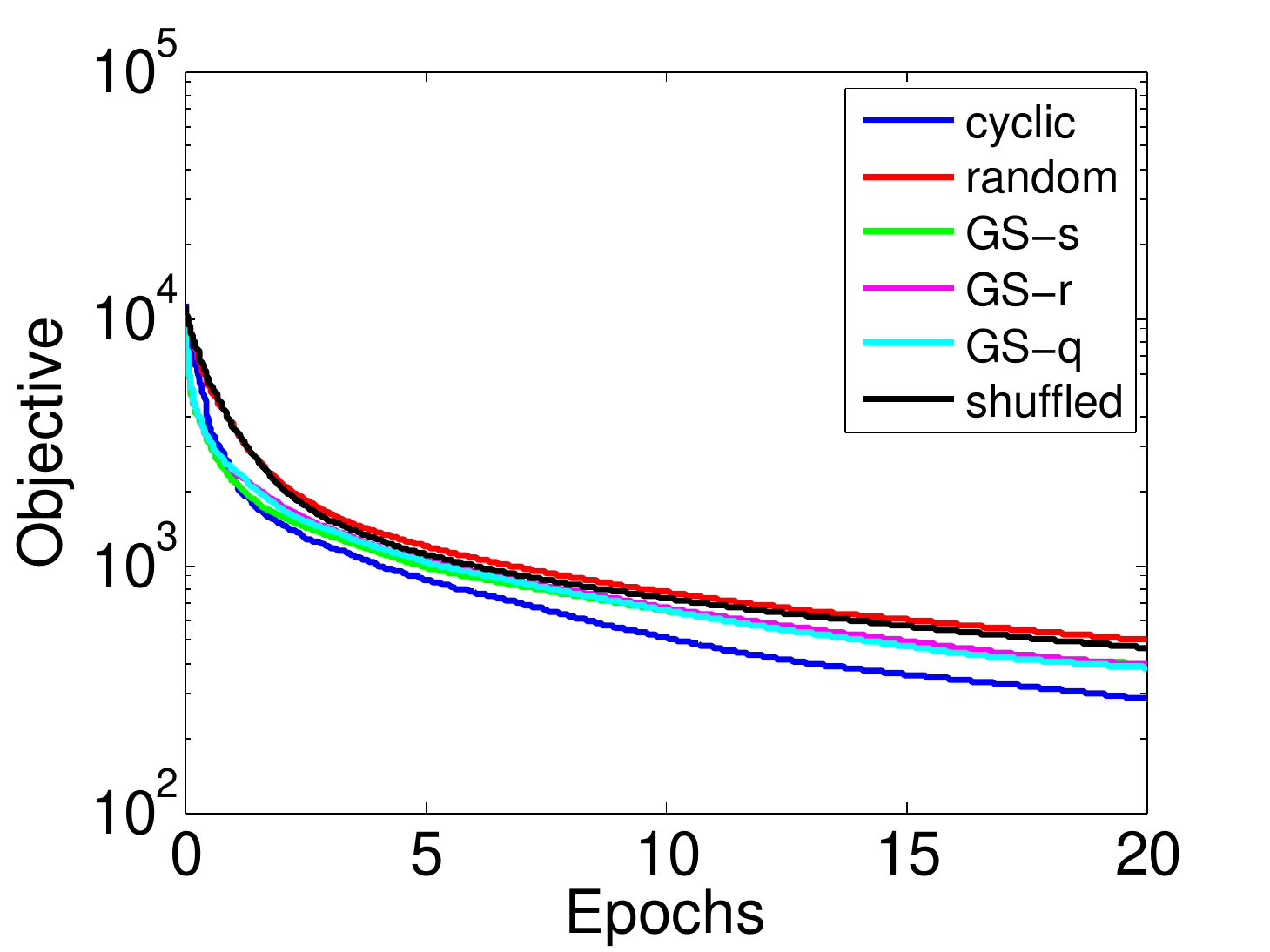}
\caption{The graph shows the objective decrease for each iteration of CD with various schemes of coordinate selection on solving the sparse logistic regression.}\label{obj_lr}
\end{figure}

We generate $m=100$ data points from two different bivariate normal distributions and apply CD to \eqref{eq:LR}. The objective decrease is plotted in Figure \ref{obj_lr}. The randomized and shuffled cyclic variants were averaged over 100 trials.

Note that cyclic CD outperforms all other variants, including all of the Gauss-Southwell rules. Peculiarly, shuffled cyclic performs significantly worse than cyclic CD, which is worthy of investigation.

\subsection{Support Vector Machines for Classification}\label{SVM}

We consider another popular classification model, support vector machines (SVM) \cite{SuykensVandewalle1999_LeastSquares}, to predict categories of new observations. 
Different from the logistic regression that estimates the probability of a binary response, the SVM is a non-probabilistic classifier, and it can be 
formulated as the following optimization problem:
\begin{equation} \label{eq:SVM}
\begin{aligned}
& \minimize \limits_{\vw,b, \xi} & & \frac{1}{2} \| \vw \|_2^2 + C \sum \limits_{i=1}^{m} \xi_i \\
& \text{subject to } & & y_i ( \vw \cdot \vx_i - b) \geq 1-\xi_i,\ \xi_i\ge0, \quad\forall i.
\end{aligned}
\end{equation}
where $\xi_i$'s are slack variables, 
and $C > 0$ is a penalty parameter. For an optimal solution $(\vw,b,\mathrm{\xi})$ to \eqref{eq:SVM}, it is easy to see that $\xi_i=\max (1-y_i (\vw^T \bx_i-b),0),\,\forall i.$ 

This problem can be interpreted as finding the best separating hyperplane that classifies the data by maximizing the margin width from the hyperplane to the nearest data point of either class.
To simplify the problem, we enforce $b=0$ that corresponds to the unbiased case, and we reformulate and solve its dual problem by substituting $\vw = \sum \limits_{i} \alpha_i y_i \bx_i$:
\begin{equation*} \label{eq:SVM_D}
\begin{aligned}
& \minimize \limits_{\bal} & & \frac{1}{2} \bal^T \vQ \bal - \mathbf{1}^T \bal \\
& \text{subject to } & & 0 \leq \bal_i \leq C,  \quad\forall i,
\end{aligned}
\end{equation*}
where $\bal \in \mathbb{R}^m$ is the dual variable, and $Q_{ij} = y_i y_j \vx_i^T \vx_j$.

\subsubsection{Update Derivation}


Let
\begin{equation*}
F(\bal) = \frac{1}{2} \bal^T \vQ \bal - \mathbf{1}^T \bal.
\end{equation*}
Then our goal is to solve 
\begin{equation}
\minimize_{\bal} F(\bal) + \iota_{[0, C]}(\bal).
\end{equation}
where $\iota_{[0, C]}(\bal)$ denotes the indicator function over the feasible set $\{\bal : 0 \leq \bal_i \leq C,\,\forall i\}$. We apply the prox-linear update \eqref{prox-linear} to solve this problem.

The prox-linear update solves the subproblem
$$\bal_{i_k}^k = \argmin_{\bal_{i_k}} \langle \vQ_{i_k, :} \bal^{k - 1} - 1, \bal_{i_k} - \bal_{i_k}^{k - 1} \rangle + \frac{L_{i_k}}{2} \|\bal_{i_k} - \bal_{i_k}^{k - 1}\|^2 + \iota_{[0, C]}(\bal_{i_k}).$$

The first-order optimality condition of this minimization is
\begin{align*}
& 0 \in \vQ_{i_k, :} \bal^{k - 1} - 1+ L_{i_k} (\bal_{i_k} - \bal_{i_k}^{k - 1}) + \partial \iota_{[0, C]}(\bal_{i_k}),\\
\mbox{equivalently,}\quad & \bal_{i_k}^{k - 1} - \frac{1}{L_{i_k}} (\vQ_{i_k, :} \bal^{k - 1} - 1) \in (I + \frac{1}{L_{i_k}} \partial \iota_{[0, C]})(\bal_{i_k}),
\end{align*}
and thus
\begin{align*}
& \vx_{i_k}^k = \prox_{\iota_{[0, C]}} (\bal_{i_k}^{k - 1} - \frac{1}{L_{i_k}} (\vQ_{i_k, :} \bal^{k - 1} - 1))
\end{align*}
as desired. Since the proximal operator of the indicator function is the projection operator, this yields the projected coordinate update

\begin{equation}\label{eq:svm-pj-cgrad}
\bal_{i_k}^{k+1} = \min(\max(\bal_{i_k}^{k}- \frac{\grad_{i_k} f(\bal^{k})}{Q_{i_k,i_k}},0),C).
\end{equation}


\subsubsection{Derivations for Gauss-Southwell Rules}
We derive the update schemes for the Gauss-Southwell rules:
\begin{enumerate}
\item \textit{GS-s index rule}: Recall that the GS-s rule chooses an index $i_k$ at each iteration $k$ by (\ref{gs-s}). Note for the SVM problem, $f(\bal) = \frac{1}{2} \bal^T \vQ \bal - 1^T \bal$, and thus
$$ \grad f(\bal^{k-1}) = \vQ \bal^{k-1} - 1.$$
In addition, for any $\bal$ with $0\le\alpha_i\le C,\,\forall i$, we have that any $\vg\in\partial \iota_{[0,C]}(\bal)$ satisfies
$$g_i\left\{\begin{array}{ll}\le 0, &\text{ if }\alpha_i=0,\\
=0, &\text{ if }0<\alpha_i<C,\\
\ge 0, &\text{ if }\alpha_i=C,\end{array}\right.$$
for all $i$.
Therefore, by \eqref{gs-s}, 
the scores are computed by
$$S_j(\bal^{k - 1}) = \begin{cases}
\min(\vQ_{j, :} \bal^{k - 1} - 1, 0), & \mbox{if $\bal^{k - 1}_j = 0$}, \\
\max(\vQ_{j, :} \bal^{k - 1} - 1, 0), & \mbox{if $\bal^{k - 1}_j = C$}, \\
 \vQ_{j, :} \bal^{k - 1} - 1, & \mbox{otherwise.} \end{cases}$$
We choose the index corresponding to the largest score, i.e., $i_k=\argmax_j|S_j(\bal^{k - 1})|$.
\item \textit{GS-r index rule}:  Recall that the GS-r rule chooses an index $i_k$ at each iteration $k$ by (\ref{gs-r}). Since the proximal operator of the indicator $\iota_{[0, C]}$ is the projection operator, the index is given by
$$i_k = \argmax_j \big|\bal_j^{k - 1} - \min(\max(\bal_j^{k - 1} - \frac{1}{\vQ_{j, j}} \grad_j f(\bal^{k - 1}),0), C)\big|. $$
\item \textit{GS-q index rule}: Recall that the GS-q rule chooses an index $i_k$ at each iteration $k$ by (\ref{gs-q}). Let
$$ d = \min(\max(\bal^{k - 1}_j - \frac{1}{\vQ_{j, j}} \grad_j f(\bal^{k - 1}),0),C)-\bal^{k-1}_j.$$
Plugging $d$ into the quadratic approximation
$$ f(\bal^{k-1}) + \nabla_j f(\bal^{k-1})^T d + \frac{L_j}{2} \|d\|^2$$
yields the greedy scores. Choosing the index corresponding to the smallest score yields the index $i_k$.
\end{enumerate}

\subsubsection{CF Analysis}

It is easy to see that one coordinate update (\ref{eq:svm-pj-cgrad}) for SVM costs $O(m)$, mainly in evaluating $\vQ_{i_k,:}\bal^k$, while a full projected gradient update needs to first compute the gradient $\vQ\bal^k$ that costs $O(m^2)$ and then project at an additional cost of $O(m)$. Hence, by the definition of CF property discussed in section \ref{cf structures}, the update (\ref{eq:svm-pj-cgrad}) is CF.

For the coordinate update, we cache $\vQ \bal^{k-1}$ and update $\bal^k$ as
\begin{equation*}
\bal_{i_k}^k = \min(\max(\bal^{k-1}_{i_k}-\frac{(\vQ \bal^{k-1})_{i_k}-1}{Q_{i_k,i_k}},0),C).
\end{equation*}
We then maintain $\vQ \bal^{k}$ by
\begin{equation*}
\vQ \bal^{k} = \vQ \bal^{k-1}+\vQ_{i_k,:}(\bal^{k}_{i_k}-\bal^{k-1}_{i_k}).
\end{equation*}
Both steps take $O(m)$ operations, so this way also shows the update (\ref{eq:svm-pj-cgrad}) is CF.

\subsubsection{Numerical Example}
\begin{figure}
\center
\includegraphics[scale=0.6]{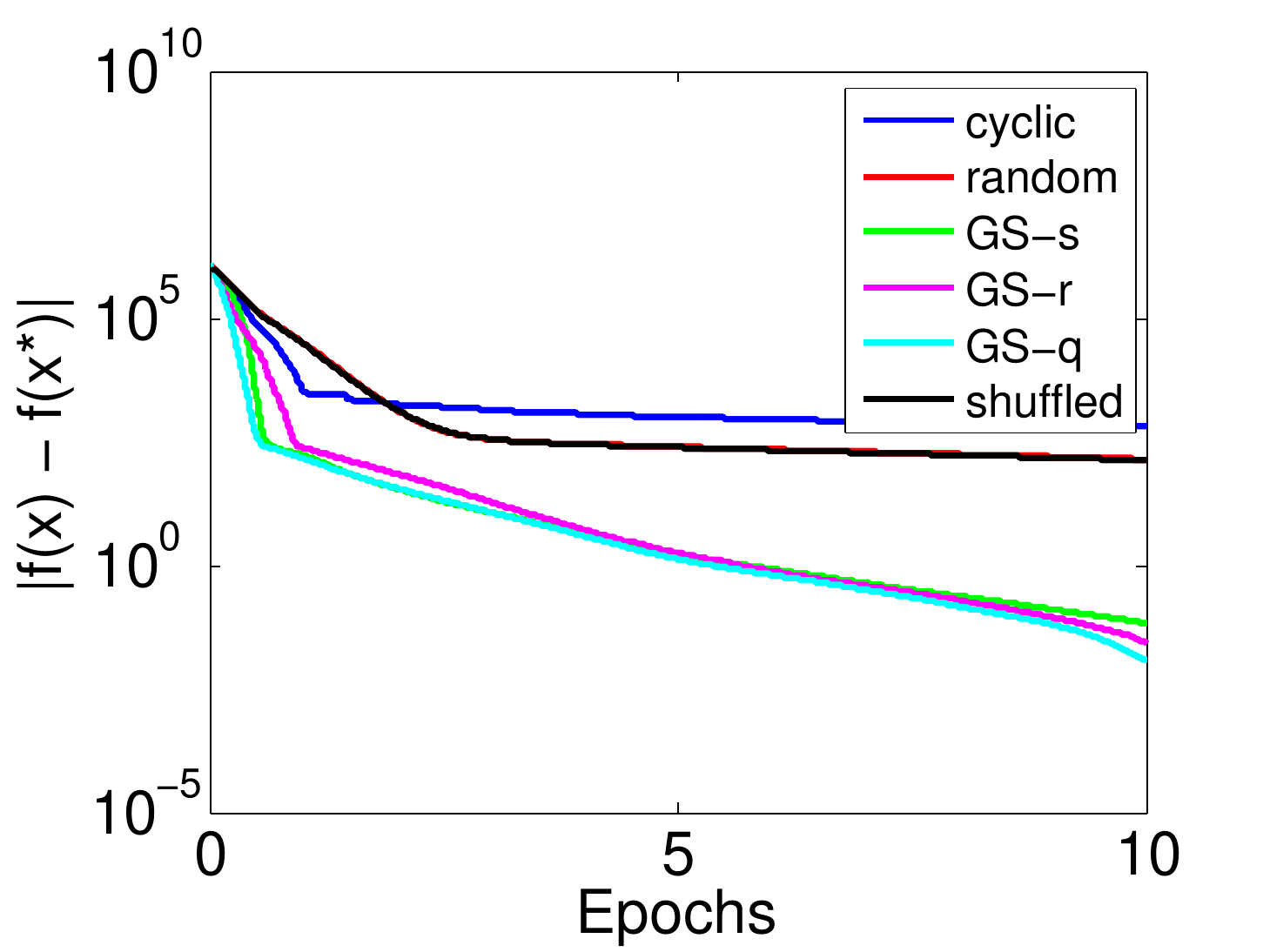}
\caption{The graph shows the objective decrease with respect to the number of iterations by applying CD to SVM on the a2a dataset.}\label{obj_svm}
\end{figure}

We train a classifier by applying our approach to the a2a training set, consisting of 20,242 training examples with 123 features \cite{chang2011libsvm,Lichman:2013}. We average our results over 10 trials. We use the update given in \eqref{eq:svm-pj-cgrad} and compare cyclic, randomized, shuffled cyclic, GS-s, GS-r, and GS-q index rules. 
The convergence results for the objective decrease over 10 epochs are plotted in Figure \ref{obj_svm}.

We note, in particular, that the GS-q and GS-s rule perform significantly better than the other rules for this application. Randomized, GS-r, and shuffled cyclic CD perform similarly.

\subsection{Semidefinite Programming}

A standard semidefinite program (SDP) is of the following form:
\begin{equation}\label{eq:sdp}
\begin{aligned}
& \underset{\vX}{\text{minimize}}
& & \langle \vC, \vX \rangle \\
& \text{subject to}
& & \mathcal{A}(\vX) = \vb, \\
& & & \vX \succeq 0,
\end{aligned}
\end{equation}
where $\cA$ is a linear operator.
Typically, SDPs can be solved in polynomial time using interior-point methods, but in practice, large-scale SDPs require an enormous amount of work at each iteration if an interior-point method is used. Because of the increasing size of the SDPs encountered in modern applications, \cite{wen2012block} develops a block coordinate descent approach that can solve large-scale SDPs much more cheaply per iteration than interior-point methods.

The procedure in \cite{wen2012block} may be interpreted as a row-by-row block minimization method. It can be summarized as cyclically updating the rows and columns of $\vX$, one pair at a time, by minimizing the objective of \eqref{eq:sdp} and keeping the constraints satisfied. Namely, let $\vX^{k, i}$ denote the current value of $\vX$ before performing the $i$th inner update in the $k$th outer iteration. At a given outer iteration $k$ and inner iteration $i$, the algorithm seeks to update $\vX$ to $\vX^{k, i+1}$ by solving the following optimization problem:
\begin{equation*}
\begin{aligned}
& \underset{\vX}{\text{minimize}}
& & \langle \vC, \vX \rangle \\
& \text{subject to}
& & \mathcal{A}(\vX) = \vb, \\
& & & \vX \succeq 0,\\
& & & \vX_{\neq i, \neq i} = \vX_{\neq i, \neq i}^{k, i},
\end{aligned}
\end{equation*}
where $\vX_{\neq i, \neq i}$ denotes the submatrix of $\vX$ excluding its $i$th row and column.  
This method, called the RBR method, has theoretical guarantees for SDPs with simple bound constraints, but may also be applied to more general bound constraints. 
We refer interested readers to \cite{wen2012block} for more information about the method.

\section{Implementations for Large-Scale Systems} \label{parallel}

Although serial CD performs well by relying on cheap updates, serial CD may not be capable of solving some large-scale systems due to memory and speed constraints in the problem. In this case, additional speedup must be gained by implementing portions of the solver in a parallel or distributed fashion. Many large-scale applications already utilize parallelisms, such as video processing, 4D-CT processing, large-scale dynamical systems, systems with streaming data, and tensor factorizations.

Parallel computing breaks a problem into simpler parts that are executed simultaneously by multiple agents while being coordinated by a controller. By using multiple cores, CPUs, or networked computers in parallel, we can overcome potential memory and speed constraints in solving our problem.

In this section, we motivate the use of parallel and distributed computing for scaling coordinate descent algorithms for larger systems. We introduce some solutions for scaling CD for larger problems in both multicore and multi-machine architectures. We also discuss parallelized numerical linear algebra methods and parallelized CD methods and give relevant problem structures and resources for implementation.

\subsection{Parallelization of Coordinate Updates}

By leveraging the CF properties of a problem as we have discussed in Section \ref{cf structures}, coordinate descent algorithms are made computationally worthwhile by relying on cheap updates to iteratively solve optimization problems. However, as our problem scales in size, the amount of work for each coordinate block update increases, stalling the computation time of each update.

For example, consider the computation of the gradient for the least-squares problem, $\vA^T \vA \vx - \vA^T \vb$, which we have discussed earlier. If our data matrix $\vA$ has dimensions $1,000,000 \times 1,000,000$, then the cost of computing a coordinate update using $400$ coordinate blocks is equivalent to the cost of calculating a full gradient update for a problem with dimensions $50,000 \times 50,000$. Therefore, to solve large-scale problems quickly, we must leverage multicore systems to gain additional speedup.

\subsubsection{Parallelized Numerical Linear Algebra}

Consider the least squares problem as discussed in Section \ref{cf structures}. Recall that the block coordinate update for least squares is given by:
$$\vx^k_{i_k} = \vx^{k - 1}_{i_k} - \alpha ((\vA^T \vA)_{i_k} \vx^{k - 1}_{i_k} - (\vA^T \vb)_{i_k})$$
In order to compute a coordinate update for least squares, assuming $\vA^T \vA \in \mathbb{R}^{m \times m}$ and $\vA^T \vb \in \mathbb{R}^m$ are precomputed and we have $s$ blocks, the cost of each operation is as follows:
\begin{enumerate}
\item Matrix-Vector Multiplication $(\vA^T \vA)_{i_k} \vx^{k - 1}_{i_k}$: $O(\frac{m^2}{s})$
\item Vector Difference $(\vA^T \vA)_{i_k} \vx^{k - 1}_{i_k} - (\vA^T \vb)_{i_k}$: $O(\frac{m}{s})$
\item Scalar-Vector Multiplication $\alpha((\vA^T \vA)_{i_k} \vx^{k - 1}_{i_k} - (\vA^T \vb)_{i_k})$: $O(\frac{m}{s})$
\item Vector Difference $\vx_{i_k}^{k - 1} - \alpha((\vA^T \vA)_{i_k} \vx^{k - 1}_{i_k} - (\vA^T \vb)_{i_k})$: $O(\frac{m}{s})$
\end{enumerate}
Assuming no communication cost, we see that the major bottleneck in computing each coordinate update consists of numerical linear algebra operations. Therefore, we can improve the efficiency of our coordinate updates by parallelizing our numerical linear algebra operations.

Since writing stable, efficient parallel numerical linear algebra solvers may be difficult, we list some common libraries and packages that are useful for implementing parallelized numerical linear algebra and point the reader to additional reports and references readily available in the public domain:
\begin{itemize}
\item \textbf{BLAS}: Basic Linear Algebra Subprograms (BLAS) is the standard low-level routines for performing linear algebra operations. BLAS has been implemented in both sequential and parallel fashions and has been designed to be highly optimized for high-performance computing. BLAS is categorized into three levels: Level 1 consists of vector operations, Level 2 consists of matrix-vector operations, and Level 3 consists of matrix-matrix operations. Please refer to the BLAS website at \url{http://www.netlib.org/blas/} for more information.
\item \textbf{LAPACK}: Linear Algebra PACKage (LAPACK) provides routines for solving more complex linear algebra operations, such as solving systems of simultaneous linear equations, least-squares solutions of linear systems of equations, eigenvalue problems, and beyond. These libraries are designed to run on shared-memory parallel processors. Please refer to the LAPACK website at \url{http://www.netlib.org/lapack/} for more information.
\item \textbf{PLASMA}: Parallel Linear Algebra Software for Multicore Architectures (PLASMA) is a dense linear algebra package designed for multicore computing. PLASMA offers routines for solving linear systems of equations, least squares problems, eigenvalue problems, etc. Please refer to the PLASMA website at \url{http://icl.cs.utk.edu/plasma/} for more details.
\end{itemize}

\subsubsection{Parallelized Coordinate Descent}

Alternatively, we can also improve the timeliness of CD for large-scale problems by parallelizing the coordinate updates, and gain an approximate speedup (ideally) proportional to the number of processors. In particular, rather than performing each operation faster, each core performs a partial coordinate update, which together form a full coordinate or gradient update. In order for parallel CD to be effective, it again depends on the CF analysis of our problem. We list some fairly common problem structures which lend themselves well to parallelizing CD:

\begin{itemize}
\item \textit{Separability}: A function $F$ is separable if it can be written as
$$F(\vx) = \sum_{i = 1}^s f_i(\vx_i)$$
where each $f_i$ only depends on a non-overlapping coordinate block $\vx_i$ of $\vx = (\vx_1, \vx_2, ..., \vx_s)$.
\item \textit{Partial separability}: A function $F$ is partially separable if it can be written as
$$F(\vx) = \sum_{J \in \mathcal{J}} f_J (\vx)$$
where $\mathcal{J}$ is a finite collection of nonempty subsets of $\{1, ..., s\}$ and $f_J$ depends on blocks $\vx_i$ for $i \in J$ only. If
$$|J| \leq \omega \text{ for all } J \in \mathcal{J}$$
then we say that $f$ is partially separable with degree $\omega$.
\end{itemize}

Clearly, minimizing a separable function $F(\vx)$ is equivalent to the independent minimization of each $f_i$ over $\vx_i$, which is obviously parallelizable since the objective is minimized if each core minimizes over each or partition of functions $f_i$.

Partially separable functions are also useful since they only couple some components of $\vx$ together. Some common examples of partially separable functions include:
\begin{enumerate}
\item Square Loss: $f_j(\vx, \vA_j, \vy_j) = \frac{1}{2} (\vA_j^T \vx - \vy_j)^2$
\item Logistic Loss: $f_j(\vx, \vA_j, \vy_j) = \log (1 + e^{-\vy_j \vA_j^T \vx})$
\item Hinge Square Loss: $f_j(\vx, \vA_j, \vy_j) = \frac{1}{2} \max\{0, 1 - \vy_j \vA_j^T \vx\}^2$
\end{enumerate}
where $\vA_j \in \mathbb{R}^n$ is a training example with label $\vy_j \in \mathbb{R}$ and $F(\vx) = \sum_{j = 1}^m f_j(\vx, \vA_j, \vy_j)$. As $\vA$ is a sparse matrix, each example may depend only on a few features, and thus the objective is partially separable. The maximum number of dependencies over all examples is then the degree of partial separability $\omega$ \cite{richtarik2016parallel}.

When parallelizing CD, we typically solve problems of the form
$$\minimize_{\vx} ~ F(\vx) = f(\vx) + r(\vx)$$
where $f$ is a partially separable smooth convex function and $r$ is a simple separable convex function, also usually proximable. This is an example of a problem that parallelizes well with CD.

Many parallel CD implementations using partial coordinate updates have already been investigated with various update schemes, sampling schemes, or step sizes, exploiting certain coordinate-friendly structures, or reducing operations by maintaining additional quantities involved in the algorithm. We refer the reader to additional work done in \cite{bradley2011parallel, fercoq2013smooth, fercoq2015accelerated,jaggi2014communication, marevcek2015distributed,  necoara2013efficient, peng2013parallel,richtarik2015optimal,richtarik2016parallel,richtarik2016parallel,tappenden2015complexity} for more detail on specific implementations.

\paragraph{Synchrony and Asynchrony}

Performing partial coordinate updates with each core may also be interpreted as performing coordinate updates on a finer partition of coordinate blocks, with a synchronization step after each set of updates. The relaxation of this synchronization gives another set of methods for parallelizing CD: asynchronous parallel CD algorithms.

Synchronous and asynchronous parallelism may be summarized as:

\begin{enumerate}
\item \textit{Synchronous Parallelism}: Synchronous algorithms distribute the coordinate computation across multiple agents and regularly synchronize across all agents to ensure consistency. The synchronization step consists of sharing the results of all coordinate updates across all agents before further computation.
\item \textit{Asynchronous Parallelism}: Asynchronous algorithms weaken or eliminate consistent synchronization across agents while still partitioning computation for execution in parallel on multiple agents. Each agent may compute with the possibly stale information it has, even if the results from other agents have not been received. 
\end{enumerate}

One can easily describe this difference through an intuitive example. Suppose the lead of a project would like to divide up operations between multiple employees. One approach would be to delegate a batch of operations for a set of employees, one operation for each employee, and wait until all employees' work has been completed before moving onto the next phase of the project. This setting would correspond to synchronous parallelization. Alternatively, the project lead may delegate a new operation to each employee as each finishes their work regardless of the other employees' progress. This setting would correspond to asynchronous parallelization.

\begin{figure}
\center
\includegraphics[width=0.8\textwidth]{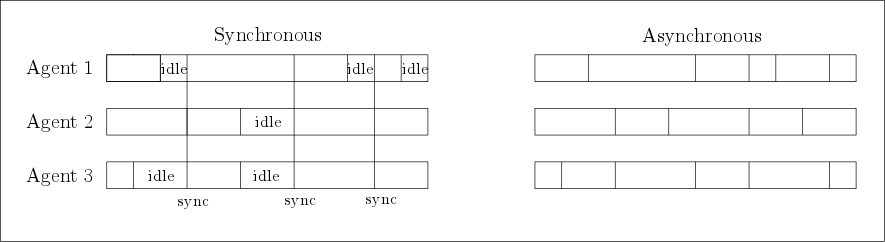}
\caption{Comparison of synchronous parallel and asynchronous parallel implementations.}\label{sync async}
\end{figure}

We highlight some major differences between synchronous and asynchronous CD methods. In particular, synchronous CD methods have been studied for much longer and are much more reliable. They have already been implemented for a variety of systems, and software is publicly available. However, synchronization requires every core, no matter how efficient, to wait for the slowest core to be communicated. Concurrent data exchanges during synchronization lead to slowdown due to lock contention and bus contention. Consequently, the speedup factor of synchronous parallel algorithms seldom reaches the number of cores or processors.

On the other hand, asynchronous methods relax the synchronization step and therefore reduce idle time and contentions that synchronization creates. In asynchronous methods, each core or processor instead performs computations using whatever information it has, regardless of the current state of the data. As a result, asynchronous methods have the potential to perform much faster - in one experiment on solving sparse logistic regression, demonstrating a 25x speedup on 32 cores rather than 4x speedup with synchronization. However, asynchrony introduces input delay or age to the components of $\vx$ used. This change makes it challenging to determine rigorous analysis. Since components are constantly updated in asynchronous computing, line search cannot be used to select a step size. Also, very few open-source software packages are available for asynchronous coordinate methods at this time. A recent open-source package is TMAC \cite{EdmundsPengYin2016_TMACToolbox}.

Figure \ref{sync async} summarizes the main difference between synchronous and asynchronous approaches \cite{wright2015coordinate}. Asynchronous CD variants differ in the assumptions they make on the choice of update components $i_k$, on the step size, and on the ``age" of the components of $\vx^k$ used. We refer the reader to additional work in \cite{bertsekas1989parallel, liu2015asynchronous, liu2015asynchronous2,HsiehYuDhillon2015_passcode,peng2015arock,HannahYin2016_UnboundedDelay} for discussion on specific asynchronous variants.

\subsubsection{Resources}

We briefly list some popular programming resources for implementing parallelized coordinate methods on a multicore machine or cluster.

\begin{itemize}
\item \textbf{Multithreading}: Multiple threads can also be executed concurrently by one or more cores while sharing memory resources. This introduces thread-level parallelism and increases utilization of the cores. The usage and syntax for threads differ for various programming languages but are supported in common languages such as C/C++, Java, and Python.
\item \textbf{OpenMP}: Open Multi-Processing (OpenMP) is an API for multi-platform shared-memory parallel programming in C/C++ and Fortran. It has its set of compiler directives, library routines, and environment variables that influence run-time behavior, and provides a simple interface for developing parallel applications. More information may be found at \url{http://openmp.org}.
\item \textbf{MPI}: Message Passing Interface (MPI) is a standardized message-passing system for parallel computing. There are several efficient open-source implementations of MPI available for parallel software development, supported in common languages like C/C++, Java, Python, Matlab, and R.
\end{itemize}

\section{Conclusion}\label{conclusion}

Coordinate descent has become an important optimization tool used to solve many problems arising from machine learning and large data analysis. We introduced and surveyed modern coordinate descent methods, including both elementary and block settings, for engineers and practitioners. We gave relevant theory and examples to help practitioners implement CD for various applications. Interested readers may refer to the bibliography for further elaborations and extensions on the topics explored in this monograph. We expect new adaptations and variants on coordinate descent methods as well as the new theory for understanding the nuances of CD to be developed as this class of algorithms become more understood and utilized in the major application areas.

\section*{Acknowledgements}

We would like to thank Charlotte Abrahamson, Yan Dong, Brent Edmunds, Xiaoyi Gu, Robert Hannah, Zhimin Peng, Tianyu Wu, and Wenli Yan for their helpful and insightful comments.

\newpage

\appendix

\section{Modeling Using Extended Valued Functions}\label{sec: extended value}

It is often useful in optimization to reformulate problems using \textit{extended valued functions}. They may be used to incorporate domain of functions or feasible set constraints in the objective function. We give a precise treatment below.

\begin{definition}\label{extended value}
An extended valued function is a function that maps to elements in the extended real line $f: \vX \mapsto \bar{\RR} = \RR\cup \{-\infty, \infty\}$.
\end{definition}

To motivate the use of extended valued functions in optimization, we present some examples.

\begin{example}\label{indicator}
Let $\mathcal{X} \subset \RR^n$. Then the indicator function of $\mathcal{X}$, defined as
$$\iota_{\mathcal{X}} (\vx) = \begin{cases}
0 & \mbox{ if } \vx \in \mathcal{X}\\
\infty & \mbox{otherwise,}
\end{cases}$$
is an extended valued function.
\end{example}

Indicator functions allow us to write constraints into the objective function and treat our problem as an unconstrained minimization problem. In particular, consider the problem
\begin{align*}
& & \minimize_{\vx} ~ & f(\vx) &\\
& & \subjectto ~ & \vx \in \mathcal{X} &
\end{align*}
where $f$ is a convex function and $\mathcal{X} \subset \RR^n$ is a convex set. Then our constrained problem may be rewritten as
$$\minimize_{\vx} ~ f(\vx) + \iota_{\mathcal{X}}(\vx).$$ 

In addition, we can use extended valued functions to ignore the domain of the function. 

\begin{example}\label{domain constraints}
Consider the problem
$$\minimize_{x > 0} ~ f(x) = \frac{1}{\sqrt{x}} + x.$$
We can define a new function
$$\tilde{f}(x) = \begin{cases}
\frac{1}{\sqrt{x}} + x & \mbox{ if } x > 0\\
\infty & \mbox{ otherwise}
\end{cases}$$
to remove implicit domain constraints from the function. Since $f$ and $\tilde{f}$ share the same set of minimizers, it is sufficient to consider the minimization of $\tilde{f}$.
\end{example}

By similarly rewriting optimization problems using extended valued functions, we can ignore constraints and instead optimize an unconstrained extended-valued problem, and therefore apply unconstrained minimization techniques, such as coordinate descent methods. 

\section{Subdifferential Calculus}\label{section: subdifferential calculus}

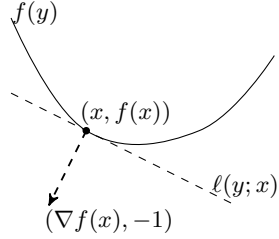
\begin{figure}
\center
	\begin{tikzpicture}[auto]
		\draw plot[smooth, tension=.7] coordinates {(-4.5,3.5) (-3.5,2) (-2,2)  (-1,3)};
		\node[smallnode]  at (-3.5,2) {};
		\node  at (-2.9757,2.2368) {{$(x,f(x))$}};
		\node  at (-4.1773,3.5703) {{$f(y)$}};
		\draw[dashed] (-4.5,2.5) -- (-1.5,1);
		\node at (-1.3811,1.2674) {{$\ell(y;x)$}};
		\draw[style=fwddash] (-3.5,2) -- (-4,1);
		\node  at (-3.2,0.8) {{$(\nabla f(x),-1)$}};
	\end{tikzpicture}
	\caption{Geometric interpretation of the definition of convexity for differentiable functions.}\label{gradient-convexity}
\end{figure}

Recall that for differentiable convex functions, we have that
$$f(\vy) \geq f(\vx) + \grad f(\vx)^T (\vy - \vx)$$
for all $\vx, \vy \in \RR^n$. In words, the line $\ell(\vy; \vx) = f(\vx) + \grad f(\vx)^T (\vy - \vx)$ is a linear underestimator of the function $f$ at $\vx$. This interpretation of convexity is shown in Figure \ref{gradient-convexity}. In particular, this inequality is equivalent to
$$\langle (\grad f(\vx), -1), (\vy - \vx, f(\vy) - f(\vx)) \rangle = \langle \grad f(\vx), \vy - \vx \rangle + f(\vx) - f(\vy) \leq 0.$$
In words, the line connecting $(\vx, f(\vx))$ and $(\vy, f(\vy))$ makes an obtuse angle with the vector $(\grad f(\vx), -1)$, as shown in the figure.

This definition of convexity motivates a more general notion of the gradient that applies to non-differentiable convex functions, called a \textit{subgradient}. 

\begin{definition}
A subgradient at $\vx \in \dom f$ is any element $\vg \in \RR^n$ such that
$$f(\vy) \geq f(\vx) + \vg^T (\vy - \vx), ~ \forall \vy \in \dom f.$$
The subdifferential $\partial f(\vx)$ is the set of all subgradients at $\vx$, i.e.
$$\partial f(\vx) := \{\vg \in \RR^n : f(\vy) \geq f(\vx) + \vg^T (\vy - \vx), ~ \forall \vy \in \dom f\}.$$
\end{definition}

Note that subdifferentials are nonempty for proper convex functions in the interiors of their domains. Subdifferentials may be empty on the domain boundaries for convex functions and anywhere for non-convex functions.

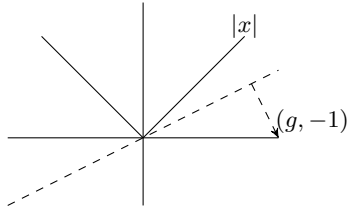
\begin{figure}
\center
	\begin{tikzpicture}[scale = .45]
		\draw [thin] (-4,0) -- (4,0);
		\draw [thin] (0,-2) -- (0,4);
		\draw [thin] (-3,3) -- (0,0) -- (3,3);
		\node at (3, 3.3) {{$|x|$}};
		\draw [dashed] (-4,-2) -- (4,2);
		\draw [dashed][->] (3.2,1.6) -- (4,0);
		\node at (5, 0.5) {{$(g,-1)$}};
	\end{tikzpicture}
\caption{A subgradient of $f(x) = |x|$ at $x = 0$.}\label{subgradient}
\end{figure}

Subgradients and subdifferentials help guide the development and analysis of optimization algorithms for non-differentiable functions. The most common example of a non-differentiable convex function is the absolute value function $|x|$. The subdifferential of $|x|$ is
$$\partial |x| = \begin{cases}
\{ 1 \} & \mbox{ if } x > 0\\
[-1, 1] & \mbox{ if } x = 0\\
\{ -1 \} & \mbox{ if } x < 0
\end{cases}.$$
A subgradient of $|x|$ at $x = 0$ is shown in Figure \ref{subgradient}.

Another common example is the indicator function of a closed nonempty convex set $C$. Then by definition, the subdifferential of the indicator function at $\vx \in C$ is
\begin{align*}
\partial \iota_{C}(\vx) & = \{\vg \in \RR^n : \iota_{C}(\vy) \geq \iota_C (\vx) + \vg^T(\vy - \vx), ~ \forall \vy \in \RR^n\}\\
& = \{\vg \in \RR^n : \vg^T(\vy - \vx) \leq 0, ~ \forall \vy \in C\}\\
& := N_C(\vx),
\end{align*}
also called the \textit{normal cone}. Note if $\vx \notin C$, then the subdifferential $\partial \iota_C(\vx) = \emptyset$.

We state some well-known properties of subgradients and subdifferentials.

\begin{properties*}

Assume all functions are proper convex, and $\phi_i$ is differentiable for all $i$.

\begin{enumerate}
\item If $f$ is differentiable, then $\partial f(\vx) = \{\grad f(\vx)\}$.
\item If $f(\vx) = c_1 f_1 (\vx) + c_2 f_2(\vx)$ with $c_1, c_2 \geq 0$, then
$$\partial f(\vx) \supseteq c_1 \partial f_1 (\vx) + c_2 \partial f_2 (\vx).$$
\item If $f(\vx) = h(\vA \vx + \vb)$, then
$$\partial f(\vx) \supseteq \vA^T \partial h(\vA \vx + \vb).$$
\item If $\lambda \geq 0$ and $f(\vx) = h(\lambda \vx)$, then
$$\partial f(\vx) = \lambda \partial h(\lambda \vx).$$
\item If $f(\vx) = \max \{\phi_1 (\vx), ..., \phi_m (\vx)\}$, then for $I(\vx) = \{i : \phi_i(\vx) = f(\vx)\}$,
$$\partial f(\vx) = \conv \{\grad \phi_i (\vx) : i \in I(\vx)\}$$
where $\conv\{\cdot\}$ denotes the convex hull.
\end{enumerate}

\end{properties*}

Under more technical assumptions, such as constraint qualification, Properties 2 and 3 hold with equality. For most cases in this monograph, they indeed do. We refer the reader to \cite{rockafellar2015convex} for more detail since they lie outside of the scope of this monograph.

Using subdifferentials, we can immediately characterize optimal solutions.

\begin{theorem}
If $f: \RR^n \rightarrow (-\infty, \infty]$ is proper convex, then $\vx^*$ is a global minimizer if and only if $0 \in \partial f(\vx^*)$.
\end{theorem}

\begin{proof}
By definition, $\vx^*$ is a global minimizer of $f$ if and only if for all $\vx \in \RR^n$,
\begin{align*}
& f(\vx) \geq f(\vx^*)\\
\iff & f(\vx) \geq f(\vx^*) + 0^T(\vx - \vx^*)\\
\iff & 0 \in \partial f(\vx^*).
\end{align*}
\end{proof}

\begin{theorem}
Let $f: \RR^n \rightarrow (-\infty, \infty]$ be a proper closed convex function and $C$ be a nonempty closed convex set. Then $x^* = \argmin\{f(\vx) : \vx \in C\}$ if there exists $\vp \in \partial f(\vx^*) \cap (-N_C(\vx^*))$.
\end{theorem}

\begin{proof}
By the first order optimality condition, $\vx^*$ is a global minimizer of $f(\vx) + \iota_C(\vx)$ if and only if
$$0 \in \partial f(\vx^*) + \partial \iota_C(\vx^*) = \partial f(\vx^*) + N_C(\vx^*)$$
which holds if there exists a $\vp$ such that $\vp \in -N_C(\vx^*)$ and $\vp \in \partial f(\vx^*)$, as desired.

\end{proof}

Note that unlike gradients, which can be formed by partial derivatives, partial subgradients do not necessarily form a gradient. In particular, if
$$\vp_1 \in \partial_1 f(\vx_1, ..., \vx_n), ..., \vp_n \in \partial_n f(\vx_1, ..., \vx_n),$$
then $(\vp_1, ..., \vp_n)$ may \textit{not} be in $\partial f(\vx_1, ..., \vx_n)$.

\section{Proximal Operators}\label{proximal operators}

The \textit{proximal mapping} or \textit{proximal operator} appears in many algorithms for minimizing convex, non-smooth functions. It involves a smaller minimization problem that may be solved cheaply in certain cases. 

\begin{definition} 
Given a closed, proper, and convex function $f$, the \textit{proximal operator} for $\alpha f$ is defined as
$$\prox_{\alpha f} (\vy) = \argmin_{\vx} f(\vx) + \frac{1}{2 \alpha} \| \vx - \vy \|_2^2.$$
\end{definition}

We first show that the proximal operator is a generalization of projections.

\begin{example}
Given a convex set $\mathcal{X}$, we will show that the proximal operator for the indicator variable $\iota_{\mathcal{X}} (\vx)$ is the projection onto $\mathcal{X}$.

Note that by definition,
$\prox_{\iota_\mathcal{X}} (\vx) = \argmin_{\vu} \iota_{\mathcal{X}} (\vu) + \frac{1}{2}\|\vu - \vx\|_2^2$, which is equivalent to the problem
\begin{equation*}
\begin{aligned}
& \minimize_{\vu} & & \frac{1}{2}\|\vu - \vx\|_2^2\\
&\subjectto & & \vu \in \mathcal{X}
\end{aligned}
\end{equation*}
This problem, by definition, is the projection operator, $\proj_{\mathcal{X}}(\vx)$.
\end{example}

The second prime example of proximal operators is for minimizing the $\ell_1$-norm. In particular, we can recover the shrinkage operator, as in \cite{beck2009fast}.

\begin{example}
We will show that the proximal operator for the $\ell_1$ norm is the shrinkage operator, i.e.
$$\prox_{\mu \|\cdot\|_1}(\vx)_i = \begin{cases}
x_i - \mu & \mbox{ if } x_i > \mu\\
0 & \mbox{ if } x_i \in [-\mu, \mu]\\
x_i + \mu & \mbox{ if } x_i < -\mu.
\end{cases}$$

First note that since the minimization is separable: 
$$\|\vu \|_1 + \frac{1}{2 \mu}\|\vx - \vu\|_2^2 = \sum_{i = 1}^n |u_i| + \frac{1}{2 \mu}(x_i - u_i)^2,$$ 
it is sufficient to consider 
$$\prox_{\mu |\cdot|}(x) = \argmin_u |u| + \frac{1}{2\mu}(u - x)^2.$$
By the first order optimality condition of the minimization problem,
$$0 \in \partial |u| + \frac{1}{\mu}(u - x).$$
This gives three cases:
\begin{enumerate}
\item $u > 0 \iff 0 = \mu + (u - x) \iff u = x - \mu \iff x > \mu$
\item $u = 0 \iff 0 \in \mu [-1, 1] + (0 - x) \iff x \in [-\mu, \mu]$
\item $u < 0 \iff 0 = -\mu + (u - x) \iff u = x + \mu \iff x < -\mu$.
\end{enumerate}
Combining these three cases gives the result.
\end{example}

For the sake of completeness, we list some common properties and interpretations of proximal operators that appear in the monograph. Please refer to \cite{rockafellar2015convex} for more details on proximal operators.

\begin{theorem}
Let $(I + \alpha \partial f)$ denote the operator that takes $(I + \alpha \partial f)(\vx) = \vx + \alpha \partial f(\vx)$. Then $(I + \alpha \partial f(\vx))^{-1}(\vx) = \prox_{\alpha f}(\vx)$.
\end{theorem}

\begin{proof}
Suppose $\vu = \prox_{\alpha f}(\vx)$. Then 
\begin{align*}
	\vx = \prox_{\alpha f}(\vy) & \iff 0 \in \partial f(\vx) + \frac{1}{\alpha}(\vx - \vy)\\
	& \iff 0 \in \alpha \partial f(\vx) + (\vx - \vy)\\
	& \iff \vy \in x + \alpha \partial f(\vx)\\
	& \iff \vy \in (I + \alpha \partial f)(\vx)\\
	& \iff \vx \in (I + \alpha \partial f)^{-1}(\vy).
\end{align*}
\end{proof}

\begin{theorem}
For a separable function $f(\vx, \vy) = g(\vx) + h(\vy)$, then
$$\prox_f (\vx, \vy) = (\prox_f (\vx), \prox_g (\vy)).$$
\end{theorem}

\begin{proof}
This follows from the definition of the proximal operator and that the minimization of $f$ is equivalent to minimizing $g$ and $h$ independently.

\end{proof}

\section{Proofs for Summative Proximable Functions}\label{proofs}

We give relevant propositions and proofs for our results for summative proximable functions. Recall that we are interested in evaluating the proximal operator of the form
$$\prox_{\alpha (f+g)} (\vy)=\argmin_{\vx} f(\vx) + g(\vx)+\frac{1}{2\alpha}\|\vx - \vy\|_2^2$$
where both $f$ and $g$ have inexpensive proximal operators.

\subsection{Proof for $\ell_2$-Regularized Proximable Functions}

\begin{proposition}\label{prop:proximable1}
Let $\vx \in \RR^n$. If $f(\vx)$ is a convex, homogeneous function of order 1 (i.e., $f(\alpha \vx) = \alpha f(\vx)$ for $\alpha\ge 0$) and $g(\vx):=\beta\|\vx\|_2$, then $\prox_{f+g}=\prox_{g}\circ\prox_{f}.$
\end{proposition}
\begin{lemma}\label{lm:homo} If $f(\vx)$ is a convex, homogeneous function of order 1 (i.e., $f(\alpha \vx) = \alpha f(\vx)$ for $\alpha\ge 0$), then the following results hold
\begin{enumerate}
  \item $\partial f(\vx)=\partial f(\alpha \vx)$ for any $\alpha>0$ and $\vx\in\RR^n$;
  \item $\partial f(\vx)\subset \partial f(0)$ for any $\vx\in\RR^n$.
\end{enumerate}
\end{lemma}
\begin{proof} Part 1. Let $\alpha >0$. By the chain rule, $\partial_{\vx} f (\alpha \vx)=\alpha \partial f(\alpha \vx)$. By $f(\alpha \vx) = \alpha f(\vx)$, we also have $\partial_{\vx} f (\alpha \vx) = \partial_{\vx} (\alpha f(\vx)) = \alpha \partial_{\vx} f(\vx)=\alpha \partial f(\vx)$. Hence, $\partial f(\vx)=\partial f(\alpha \vx)$.


Part 2. For any $\vx$ and $\alpha >0$, let $\vy = \alpha \vx$. Let $\vp \in \partial f(\vx)$. By part 1, $\vp\in\partial f(\vy)$. Thus, by definition, $\vp$ obeys $f(\vz)\ge f(\vy) + \langle \vp,\vz-\vy \rangle$ for any $\vz\in\RR^n$,  and this inequality holds for any $\alpha >0$. Now let $\alpha \to 0$ and, by continuity, we have $f(\vz)\ge f(0) + \langle \vp,\vz-0 \rangle$. Hence, $\vp\in\partial f(0)$, and $\partial f(\vx)\subset \partial f(0)$.
\end{proof}

\begin{proof}[Proof of Prop. \ref{prop:proximable1}]
Let $\vx \in \RR^n$, $\vy:=\prox_f(\vx)$, and $\vz:=\prox_g(\vy)$. We shall show that $\vz=\prox_{f+g}(\vx)$.

From $\vy:=\prox_f(\vx)$ and $\vz:=\prox_g(\vy)$, we obtain the optimality conditions of their minimization problems, respectively,
\begin{align*}
  0  & \in \partial f(\vy)+(\vy-\vx), \\
  0  & \in \partial g(\vz)+(\vz-\vy),
\end{align*}
and adding them gives us
\begin{align*}
  0  & \in \partial f(\vy)+\partial g(\vz)+(\vz-\vx).
\end{align*}
Now using $g(\cdot):=\beta\|\cdot\|_2$, we have $\vz:=\lambda_{\vy} \vy$, where $\lambda_{\vy} = \frac{\vy}{\|\vy\|_2}\max\{0,\|\vy\|_2-\beta\}\ge 0$. By Lemma \ref{lm:homo}, Part 1 for the case $\lambda>0$ and Part 2 for the case $\lambda=0$, we arrive at $\partial f(\vy)\subseteq \partial f(\vz)$ and thus
\begin{align*}
  0  & \in \partial f(\vz)+\partial g(\vz)+(\vz-\vx),
\end{align*}
which is the optimality condition for $\vz=\prox_{f+g}(\vx)$.
\end{proof}
In the proof, the formula of $\lambda_{\vy}$ is not important; only $\lambda_{\vy}\ge 0$ is. Therefore, Prop. \ref{prop:proximable1} remains valid for $g(\vx):=\|\vx\|_2^p$ for any $p\ge 1$.

\subsection{Proof for TV-Regularized Proximable Functions}
\begin{proposition}\label{prop:proximable2}
Let $\vx \in \mathbb{R}^n$. Define the total variation semi-norm: $\mathrm{TV}(\vx) := \sum_{i = 1}^{n - 1} |x_{i + 1} - x_i|$. If $f(\vx) = \beta \mathrm{TV}(\vx)$ and $g(\vx)$ is a closed, proper, convex function that satisfies
\begin{subequations}\label{eq:monoprox}
\begin{align}
x_i > x_{i+1} & \implies \prox_g(\vx)_i \ge \prox_g (\vx)_{i+1}\\
x_i < x_{i+1} & \implies \prox_g(\vx)_i \le \prox_g (\vx)_{i+1}\\
x_i = x_{i+1} & \implies \prox_g(\vx)_i = \prox_g (\vx)_{i+1},
\end{align}
\end{subequations}
then $\prox_{f + g} = \prox_{f} \circ \prox_g$.
\end{proposition}
\begin{proof}
Let $\vx \in \RR^n$, $\vy:=\prox_f(\vx)$, and $\vz:=\prox_g(\vy)$. We shall show that $\vz=\prox_{f+g}(\vx)$.

From $\vy:=\prox_f(\vx)$ and $\vz:=\prox_g(\vy)$, we obtain the optimality conditions of their minimization problems, respectively,
\begin{align*}
  0  & \in \partial f(\vy)+(\vy-\vx), \\
  0  & \in \partial g(\vz)+(\vz-\vy),
\end{align*}
and adding them gives us
\begin{align}\label{eq:prop2half}
  0  & \in \partial f(\vy)+\partial g(\vz)+(\vz-\vx).
\end{align}

By definition, $f(\vy) = \sum_{i = 1}^{n - 1} |y_{i + 1} - y_i|$, which satisfies $\partial f(\vy) = \sum_{i=1}^{n-1} \partial_\vy |y_{i + 1} - y_i|$. In this Minkovski sum, each term satisfies
\begin{align}\label{eq:3cases}
  \partial_{(y_i,y_{i+1})} |y_{i + 1} - y_i| =
  \begin{cases}
    (1,-1), & \mbox{if } y_i > y_{i+1} \\
    (-1,1), & \mbox{if } y_i < y_{i+1} \\
    \{\alpha (1,-1)+(1-\alpha)(-1,1):\alpha\in[0,1]\} & \mbox{otherwise},
   \end{cases}
\end{align}
where the ``otherwise'' case is the convex hull of the first two cases. Since $g$ satisfies \eqref{eq:monoprox} and $\vz=\prox_g(\vy)$, from $(y_i,y_{i+1})$ to $(z_i,z_{i+1})$ it holds that
\begin{align*}
y_i > y_{i+1} & \implies z_i \ge z_{i+1}\\
y_i < y_{i+1} & \implies z_i \le z_{i+1}\\
y_i = y_{i+1} & \implies z_i = z_{i+1}.
\end{align*}
Hence, either $(z_i,z_{i+1})$ holds for the same case of $(y_i,y_{i+1})$ in \eqref{eq:3cases}, or it belongs to the ``otherwise" case, which is a superset of the first two cases. Therefore, $\partial_{(y_i,y_{i+1})} |y_{i + 1} - y_i| \subseteq \partial_{(z_i,z_{i+1})} |z_{i + 1} - z_i|$ and thus $\partial f(\vy)\subseteq \partial f(\vz)$. This together with \eqref{eq:prop2half} yields
$$ 0 \in \partial f(\vz) + \partial g(\vz) + (\vz - \vx),$$
which is the optimality condition for $\vz = \prox_{f + g}(\vx)$.
\end{proof}

\newpage

\bibliographystyle{abbrvnat}
\bibliography{bibl,wotaoyin}

\end{document}

%% file: macros.tex
\newcommand{\vb}{{\mathbf{b}}}

\newcommand{\vd}{{\mathbf{d}}}

\newcommand{\vg}{{\mathbf{g}}}

\newcommand{\vp}{{\mathbf{p}}}

\newcommand{\vu}{{\mathbf{u}}}

\newcommand{\vw}{{\mathbf{w}}}
\newcommand{\vx}{{\mathbf{x}}}
\newcommand{\vy}{{\mathbf{y}}}
\newcommand{\vz}{{\mathbf{z}}}

\newcommand{\vA}{{\mathbf{A}}}

\newcommand{\vC}{{\mathbf{C}}}
\newcommand{\vD}{{\mathbf{D}}}

\newcommand{\vG}{{\mathbf{G}}}

\newcommand{\vM}{{\mathbf{M}}}

\newcommand{\vQ}{{\mathbf{Q}}}

\newcommand{\vX}{{\mathbf{X}}}
\newcommand{\vY}{{\mathbf{Y}}}

\newcommand{\cA}{{\mathcal{A}}}

\newcommand{\RR}{\mathbb{R}}

\newcommand{\sign}{\mathrm{sign}}

\newcommand{\supp}{{\mathrm{supp}}} 
\newcommand{\dom}{{\mathrm{dom}}} 
\newcommand{\grad}{{\nabla}}    

\DeclareMathOperator{\shrink}{shrink} 
\DeclareMathOperator*{\argmin}{arg\,min}
\DeclareMathOperator*{\argmax}{arg\,max}

\newcommand{\bc}{\begin{center}}
\newcommand{\ec}{\end{center}}

\newcommand{\bdm}{\begin{displaymath}}
\newcommand{\edm}{\end{displaymath}}

\newcommand{\beq}{\begin{equation}}
\newcommand{\eeq}{\end{equation}}

\newcommand{\bfl}{\begin{flushleft}}
\newcommand{\efl}{\end{flushleft}}

\newcommand{\bt}{\begin{tabbing}}
\newcommand{\et}{\end{tabbing}}

\newcommand{\beqn}{\begin{align}}
\newcommand{\eeqn}{\end{align}}

\newcommand{\beqs}{\begin{align*}} 
\newcommand{\eeqs}{\end{align*}}  

\newtheorem{lemma}{Lemma}
\newtheorem{proposition}{Proposition}

%% file: tikz_setting.tex
\usetikzlibrary{fit, arrows, shapes, positioning, shadows, matrix, calc}

\tikzstyle std=[line width=0.7pt]   
\tikzstyle stdthin=[line width=0.3pt]   
\tikzstyle stdthick=[line width=1.0pt]   

\tikzstyle fwd=[line width=0.7pt, ->]   
\tikzstyle fwdthin=[line width=0.3pt, ->]   
\tikzstyle fwdthick=[line width=1.0pt, ->]   
\tikzstyle fwddash=[line width=0.7pt, dashed, ->]   

\tikzstyle bwd=[double, line width=0.3pt, ->]  
\tikzstyle refl=[double,  dashed, line width=0.2pt, ->]    

\tikzstyle{every node}=[font=\small] 

\tikzset{
  >=stealth', 
  invisible/.style={opacity=0}, 
  alt/.code args={<#1>#2#3}{\alt<#1>{\pgfkeysalso{#2}}{\pgfkeysalso{#3}}}, 
  visible on/.style={alt=#1{}{invisible}}, 
  smallnode/.style={circle, fill=black, thick, inner sep=1pt, minimum size=1.5pt}, 
  punkt/.style={
           rectangle,
           rounded corners,
           draw=black, very thick,
           text width=5.5em,
           minimum height=2em,
           text centered},
  punkt_big/.style={
           rectangle,
           rounded corners,
           draw=black, very thick,
           text width=7em,
           minimum height=2em,
           text centered},
}